\documentclass[10pt,a4paper,twoside]{amsart}

\usepackage[utf8]{inputenc}
\usepackage[english]{babel}
\usepackage[T1]{fontenc}
\usepackage{microtype, fancyhdr, lmodern, xcolor}
\usepackage{listings}

\usepackage[a4paper, hmarginratio=1:1]{geometry}

\usepackage{kerkis}
\usepackage{charter}

\usepackage{amsfonts, amsmath, amsthm, amssymb, mathrsfs, amscd, mathtools}
\usepackage[all]{xy}
\numberwithin{equation}{section}	
\usepackage{faktor}
\allowdisplaybreaks

\usepackage{url, enumitem}
\usepackage[noadjust]{cite}

\usepackage{multirow,bigdelim, caption}
\usepackage{longtable}
\usepackage{booktabs}
\usepackage{graphicx}
\graphicspath{{Figures/}}

\usepackage{hyperref}

\theoremstyle{plain}
\newtheorem{theorem}{Theorem}[section]
\newtheorem{proposition}[theorem]{Proposition}
\newtheorem{corollary}[theorem]{Corollary}
\newtheorem{lemma}[theorem]{Lemma}

\newtheorem{remark}[theorem]{Remark}

\theoremstyle{definition}
\newtheorem{definition}[theorem]{Definition}

\newtheorem{maintheorem}{Theorem}

\newcommand\Zz{\mathbb{Z}}
\newcommand\Qq{\mathbb{Q}}
\newcommand\Rr{\mathbb{R}}
\newcommand\Cc{\mathbb{C}}

\newcommand\Mod{\mathrm{Mod}} 	
\newcommand\I{\mathcal{I}} 	
\newcommand\K{\mathcal{K}} 	

\newcommand{\aut}[1]{\ensuremath{\operatorname{\text{Aut}}\left({#1}\right)}}

\newcommand\Sp[2]{\mathrm{Sp}_{#1}(#2)}	

\begin{document}
	\title{Finiteness properties of the Torelli group of surfaces with 2 boundary components}
	
	\author[Stylianakis]{Charalampos Stylianakis}
	\address{University of the Aegean, Department of Mathematics, Karlovassi, 83200, Samos, Greece}
	\email{stylianakisy2009@gmail.com}


	
	\date{ }
	
	\begin{abstract}
		In this paper we prove that the Torelli group of a surface of genus at least 3 with 2 boundary components is finitely generated. As a consequence, we answer Putman's question on the finite generation of the stabilizer subgroup of the Torelli group of a non separating simple closed curve \cite[Question 1.1]{PM2}. Furthermore, we prove that the Johnson kernel is finitely generated if the genus of the surface is at least 5 and has 2 boundary components.
	\end{abstract}
	
	\maketitle

	\section{Introduction}
	
	Let $\Sigma$ be an orientable surface. Denote by $\Mod(\Sigma)$ the mapping class group of $\Sigma$. For the following definitions and constructions see for example Farb-Margalit's book \cite{BFM}. Suppose that $\Sigma$ has genus $g\geq 1$ with at most one boundary component. The action of $\Mod(\Sigma)$ on $\mathrm{H}_1(\Sigma;\Zz)$ induces a representation
	\[ \Mod(\Sigma) \to \aut{\mathrm{H}_1(\Sigma;\Zz)} \]
	with image $\mathrm{Sp}_{2g}(\Zz)$. The kernel of the latter homomorphism is called the \emph{Torelli group} and it is denoted by $\I(\Sigma)$. The definition of the Torelli group for surfaces $S$ with more than one boundary components depends on how $S$ embeds into surfaces $\Sigma$ with at most one boundary component. Consider an embeding $\lambda: S \to \Sigma$ which induces a homomorphism $\lambda_*:\Mod(S) \to \Mod(\Sigma)$. We define $\I(\Sigma,S):= \lambda_*^{-1}(\I(\Sigma))$. It is proved that $\I(\Sigma,S)$ depends on the number of connected components of $\Sigma \setminus S$ \cite{PM1}.\\
	
	It is known that $\Mod(\Sigma)$ admits a finite presentation (see for example Hatcher-Thurston's proof \cite{HT}); it is also known by Johnson that if $\Sigma$ has genus at least 3 and at most one boundary component, then $\I(\Sigma)$ is finitely generated \cite{J1}. Apart from a few cases, if $S$ has more than one boundary component, then it is not known whether $\I(\Sigma,S)$ is finitely generated or not. For the finitely generated cases we have the following theorem, which is a consequence of \cite[Theorem 1.2]{PM1}:
	
	\begin{theorem}[Putman]
		Suppose that $S$ is a surface of genus $g\geq 3$ with $n$ boundary components. Suppose that $\Sigma$ is a surface of genus at least 3 and one boundary components such that $S \subset \Sigma$ and $\Sigma \setminus S$ has $n$ connected components. Then $\I(\Sigma,S)$ is finitely generated.
		\label{Putman_fg_thm}
	\end{theorem}
	
	If $S$ has genus at least 3 with two boundary components, then, there are two cases. Either, $\Sigma\setminus S$ has two connected components or $\Sigma\setminus S$ has one connected component. The first case satisfies the conditions of Theorem \ref{Putman_fg_thm}. The second case is proved in the following Theorem.
	
	\begin{maintheorem}
		Let $\Sigma$ be a surface of genus at least 4 with one boundary component. If $S \subset \Sigma$ is a surface of genus at least 3 with 2 boundary components and $\Sigma \setminus S$ is connected, then $\I(\Sigma, S)$ is finitely generated.
		\label{fg_torelli}
	\end{maintheorem}
	
	Theorem \ref{fg_torelli} suggests that the Torelli group of a genus $g\geq 3$ surfaces is finitely generated for every number of boundary components. Also, the proof of Theorem \ref{fg_torelli} is a consequence of Theorem \ref{first_th} in which we propose a generating set for $\I(\Sigma, S)$. What we actually prove is the following. Let $T$ be a subset of $\Mod(S)$ that it normally generates $\I(\Sigma, S)$ and let $\Gamma$ be the group generated by $T$. Fix a generating set $X$ of $\Mod(S)$. We show that for every $t \in T$ and $g \in X$, then $g t g^{-1}, g^{-1} t g \in \Gamma$. This shows that $\Gamma$ is normal in $\Mod(S)$. In Section 2 we show that $\I(\Sigma, S)$ is also normal in $\Mod( S)$. This means that $\Gamma = \I(\Sigma, S)$. We note that this strategy was performed by Johnson to prove that $\I(S)$ is finitely generated \cite{J1}.\\
	
	\begin{figure}[h]
		\begin{center}
			\includegraphics[scale=.5]{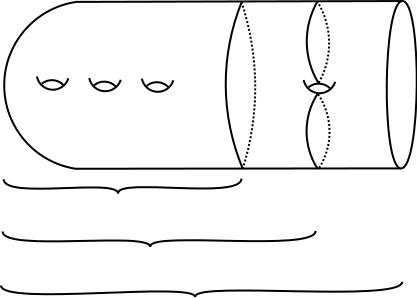}
		\end{center}
		\caption{Model of subsurfaces that is used here}
		\begin{picture}(22,12)
			\put(-40,80){$\Sigma_{g,1}$}
			\put(-20,55){$S$}
			\put(0,30){$\Sigma_{g+1,1}$}
		\end{picture}
		\label{surface_model}
	\end{figure}
	
	The generating set of Theorem \ref{first_th} is natural in the following sense. First of all, denote by $\Sigma_{g,n}$ a surface of genus $g$ with $n$ boundary components. Let $X_g$ be the Johnson's generating set of $\I(\Sigma_{g,1})$. Consider a genus $g+1$ surface $\Sigma_{g+1,1}$ with one boundary component as depicted in Figure \ref{surface_model}. From this picture, we can see that $X_g \subset X_{g+1}$. Suppose that $Y_{g+1}$ is the generating set of $\I(\Sigma_{g+1,1},S)$ Theorem \ref{first_th}. Then, $X_g \subset Y_{g+1} \subset X_{g+1}$.\\
	
	\subsection*{Stabilizer of a simple closed curve.} Theorem \ref{fg_torelli} has an immediate corollary. Let $S'$ be a surface of genus $g \geq 4$ with empty boundary and let $c$ be a nonseparating simple closed curve in $S'$. If we cut $S'$ along $c$ we end up to a surface $S_c$ of genus $g-1$ with 2 boundary components $\partial_1, \partial_2$. Denote by $\Mod(S')_c$ the stabilizer of $c$ in $\Mod( S')$. We have a short exact sequence:
	
	\[ 1 \to \langle T_{\partial_1} T^{-1}_{\partial_2} \rangle \to \Mod(S_c) \to \Mod(S')_c \to 1. \]
	
	Define the stabilizer of $c$ in $\I(S')$ as the intersection $\I(S') \cap \Mod(S')_c$ and denote it by $\I(S')_c$. By restricting the short exact sequence above to the Torelli group we have the following short exact sequence \cite[Remark in Section 2.2]{PM2}:
	
	\begin{equation}
		1 \to \langle T_{\partial_1} T^{-1}_{\partial_2} \rangle \to \I(S', S_c) \to \I(S')_c \to 1,
		\label{stabilizer_corollary}
	\end{equation}
	
	Therefore, we have the following Corollary.
	
	\begin{corollary}
		Let $c$ be a nonseparating simple closed curve of a surface $S'$ of genus $g\geq 4$. Then, $(\I(S'))_c$ is finitely generated.
	\end{corollary}
	
	\subsection*{Smaller generating sets} We now return to $\Sigma$ when there is at most one boundary component. Johnson's generating set of $\I(\Sigma)$ grows exponentially with respect to the genus $g$. Nevertheless, when $g=3$, Johnson gave the minimum number of generators; 35 for empty boundary and 42 with one boundary component. He then conjectured that for $g\geq 4$ there should be a generating set whose cardinality grows cubically with respect to the genus of $\Sigma$. That conjecture was proved by Putman \cite{PM2}, who gave a smaller generating set of size $57 \binom{g}{3}$ and that was later improved to $42 \binom{g}{3}$ by Church-Putman \cite[Proposition 4.5]{PM3}. Likewise, we also prove an analogue for a surface of genus $g\geq 3$ with 2 boundary components. Theorem B below is a consequence of Theorem \ref{handle_torelli}.
	
	\begin{maintheorem}
		Let $\Sigma$ be a surface of genus at least 4 with one boundary component. If $S \subset \Sigma$ is a surface of genus at least 3 with 2 boundary components and $\Sigma \setminus S$ is connected. Then $\I(\Sigma, S)$ is generated by $103 \binom{g}{3}$ elements.
		\label{main_cubic}
	\end{maintheorem}
	
	To prove Theorem \ref{main_cubic} we do the following. In Section \ref{small_gen_sets} we define subsurfaces $X_i \subseteq S$ with 2 boundary components and genus at most the genus of $S$. Also, $\Sigma \setminus X_i$ are connected. The subsurfaces $X_i$ are all finite in number. In Theorem \ref{cubic} we show that the union of all $\I(\Sigma,X_i)$ generate $\I(\Sigma,S)$. If the subsurfaces have genus 3, then the number of all $X_i$ is $\binom{g}{3}$. Also, in this particular case $\I(\Sigma,X_i)$ is generated by 85 elements by Theorem \ref{first_th}. Theorem \ref{main_cubic} then follows.
	
	\subsection*{Johnson kernel} Recall that $\Sigma$ is a surface of genus $g$ with at most one boundary component. Another group of interest is the following. By work of Johnson we know that the abelianization of $\I(\Sigma)$ is the direct sum of an torsion free abelian group of finite rank and a torsion vector space over $\Zz/2$ of finite rank \cite{J3}. The kernel of $\I(\Sigma) \to \mathrm{H}_1(\I(\Sigma);\mathbb{Q})$, denoted by $\mathcal{K}(\Sigma)$ is called the \emph{Johnson kernel}. Whether $\mathcal{K}(\Sigma)$ is finitely generated or not was an open problem until recently. The first positive result towards the finite generation was obtained by Dimca-Hain-Papadima \cite{DHP} who calculated the $\mathrm{H}_1(\mathcal{K}(\Sigma);\mathbb{Q})$ when the genus of $\Sigma$ is at least 6. Recently, Ershov-Sue proved that $\mathcal{K}( \Sigma)$ is finitely generated when the genus of $\Sigma$ is at least 12 \cite{ES}. The result was improved by Church-Ershov-Putman \cite{CEP} by showing that  $\mathcal{K}(\Sigma)$ is finitely generated when the genus of $\Sigma$ is at least 4. Even more recently another result was announced by Gaifullin who proved that the abelianization of $\mathcal{K}( \Sigma)$ has finite rank when the genus of $S$ is at least 3 \cite{Gaifullin}.\\
	
	Let $S\subset \Sigma$ where $S$ has two boundary components such that $\Sigma \setminus S$ is connected. We denote by $\mathcal{K}(\Sigma,S)$ the intersection $\I(\Sigma,S) \cap \mathcal{K}( \Sigma)$. Combining the work of Church \cite[Theorem 3]{TC} and Putman \cite[Theorem A]{PM4} we have that $\I(\Sigma,S) / \mathcal{K}(\Sigma,S) = \mathrm{H}_1( \I(\Sigma,S) ;  \mathbb{Q} )$. If we combine the latter result with the full abelianization of $\I(\Sigma,S)$ (which can be found in Section 2.4), then $[\I(\Sigma,S),\I(\Sigma,S)]$ has finite index in $\mathcal{K}(\Sigma,S)$.
	
	\begin{maintheorem}
		Let $\Sigma$ be a surface of genus at least 5 with one boundary component. If $S \subset \Sigma$ is a surface of genus at least 4 with 2 boundary components and $\Sigma \setminus S$ is connected, then $\mathcal{K}(\Sigma,S)$ is finitely generated.
	\end{maintheorem}
	
	In \cite{CEP} Church-Ershov-Putman proved that the Johnson kernel is finitely generated when the surface has genus at least 4 and at most one boundary component. Thus Theorem C is actually an extension of the Church-Ershov-Putman's result. The tool the authors of \cite{CEP} used to prove their Theorem is the BNS invariant. To prove Theorem C here we implicitly use the BNS invariant. Actually, we do not reference the BNS invariant here at all since our proof is explicitly based on \cite[Theorem 3.7]{CEP}. According to \cite[Theorem 3.7]{CEP} if we have groups $G < \Gamma$ such that particular conditions hold for $G$ and $\Gamma$, then $[G,G]$ is finitely generated (see Theorem \ref{CEP_THM} below). We prove that these conditions hold for $\I(\Sigma,S)<\Mod(S)$, hence $[\I(\Sigma,S),\I(\Sigma,S)]$ is finitely generated. Since the latter group has finite index in $\mathcal{K}(\Sigma,S)$, the result follows.
	
	\subsection*{Notation}
		A surface of genus $g$ with $n$ boundary components is denoted by $\Sigma_{g,n}$. In Figure \ref{surface_model} the surface $\Sigma_{g+1,1} \setminus S$ is a sphere with 3 boundary components. This latter surface is usually referred as a \emph{pair of pants} in the article. By abuse of notation when we write $f \in \Mod(S)$ we mean a representative of the mapping class $f$. Also, we write $f*g$ for $fgf^{-1}$ and $[f,g] = f g f^{-1} g^{-1}$. Finally, all modules are over $\Zz$, unless otherwise stated.

	\subsection*{Outline of the paper}
	In Section 2 we recall basics on symplectic structures and the Torelli group. In Section 3 we analyse generators of the Torelli group and we give relations between them. We combine this information to prove Theorem A. In Section 4 we prove Theorem B. In Section 5 we prove Theorem C.

	\subsection*{Acknowledgments}
	I would like to thank Andrew Putman for his directions in calculating the rational homology of the Torelli group of partitioned surfaces. Also, I would like to thank Tara Brendle and Dan Margalit for valuable suggestions and also Mikhail Ershov for suggestions on an earlier version of the article and also for giving me directions on how to improve the result of Theorem C. Finally, I would like to thank the anonymous referee from a prior submission for valuable remarks that led to substantial improvements in this work.
	
	\section{Symplectic structures and the Torelli group}
	
	The aim of this section is to recall basic definitions of symplectic structures. Then we define the Torelli group of surfaces with non empty boundary and we provide generators. In the end of the section we give the abelianization of the Torelli group.
	
	\subsection{Symplectic structures}
	
	For details on the definitions below there are a lot references in the literature, e.g. see O'Meara's book \cite{OM}. Let $R$ be a ring like $\Zz, \Rr, \Cc, \Qq, \Zz/m$. Suppose that $V$ is a module over $R$. A pairing $\iota$ is a map:
	\[ \iota ( \,, \,) : V \times V \to R \]
	which is called
	\begin{itemize}
		\item \emph{Alternating} if $\iota(x,x)=0$ for $x \in V$.
		\item \emph{Bilinear} if for $x,y,z \in V$ and $\lambda \in R$, then,
		\subitem $\iota(x+y,z) = \iota(x,z)+\iota(y,z)$,
		\subitem $\iota(x,z+y) = \iota(x,z)+\iota(x,y)$,
		\subitem $\iota(\lambda x, y) = \iota( x, \lambda y) = \lambda \iota(x,y)$.
		\item \emph{Degenerate} if there exists nonzero $x \in V$ such that $\iota(x,y) = 0$ for all $y \in V$.
	\end{itemize}
	
	Let $V$ be a module endowed with a bilinear alternating form $\iota$. Define the \emph{radical} of $V$ as $\mathrm{rad}(V) = \{ x \in V \mid \iota(x,y)=0, \, \forall \, y \in V \}$. Let $a \in V$, then $a^{\bot} = \{ x \in V \mid \iota(x,a) = 0 \}$ is \emph{orthogonal complement} with respect to $a$. If $V$ is a bilinear alternating module, then $\bar{V} = V / \mathrm{rad}(V)$ or $V = \bar{V} \oplus \mathrm{rad}(V)$. As a result, $\bar{V}$ is non-degenerate alternating module.\\
	
	A bilinear alternating non-degenerate module pairing is called a \emph{symplectic pairing}. A module $V$ endowed with a symplectic pairing is called a symplectic module. A consequence of the definitions of a symplectic module is that it has a symplectic basis $\{ x_1, \cdots, x_g, y_1, \cdots, y_g  \}$ such that
	
	\[ \iota(x_i,y_j) =
	\begin{matrix}
		0, \, \text{if } i \neq j,\\
		1, \, \text{if } i = j.
	\end{matrix}
	\]
	and $\iota(x_i, x_j) = \iota(y_i, y_j) = 0$ for any $i,j \leq g$. If $\mathrm{rad}(V)\neq 0$, then $V$ is not symplectic, and $\mathrm{rad}(V)$ has basis $\{ z_1, \cdots, z_n \}$.\\
	
	Denote by $\mathrm{Aut}(V)$ the group of automorphisms of $V$, and if $W$ is a submodule of $V$ define $\mathrm{Aut}(V,W) = \{ f \in \mathrm{Aut}(V) \mid f(W)=W \}$. Let $\bar{V}$ be a symplectic module with pairing $\iota$. Then $\mathrm{Sp}( \bar{V} ) = \{ f \in \mathrm{Aut}(\bar{V}) \mid \iota(f(x), f(y)) = \iota(x,y) \}$. Consider the following transformation of $\mathrm{Sp}( \bar{V} )$. for $a,b \in \bar{V}$, then 
	\[ a \mapsto a + \iota(a,b) b \]
	is called a \emph{transvection}. It turns out that $\mathrm{Sp}( \bar{V} )$ is generated by transvections. \\
	
	The following result is important for the proof of Theorem C. Suppose now that $\bar{V}$ is a symplectic module with a symplectic basis $\{ x_1,\cdots,x_g,y_1, \cdots, y_g \}$. We want to discuss a subgroup of $\mathrm{Sp}(\bar{V})$ stabilizing a vector in $\bar{V}$. The symplectic form induces a matrix
	\[ J = \begin{pmatrix}
		\begin{matrix} 0 & 1\\ -1 & 0 \end{matrix} & & 0 \\
		& \ddots & \\
		0 & & \begin{matrix} 0 & 1\\ -1 & 0 \end{matrix}
	\end{pmatrix} \]
	
	Then, we also have
	
	\[ \mathrm{Sp}(\bar{V}) = \{ B \in \mathrm{GL}(\bar{V}) \mid B^T J B = J \} < \mathrm{GL}(\bar{V}). \]
	
	Let $\delta$ be a vector in $\bar{V}$. Define
	
	\[ (\mathrm{Sp}(\bar{V}))_\delta = \{ B \in \mathrm{GL}(\bar{V}) \mid B^T J B = J, B\delta = \delta \} < \mathrm{Sp}(\bar{V}). \]
	
	We want to understand the structure of $(\mathrm{Sp}(\bar{V}))_\delta$. Since $\mathrm{Sp}(\bar{V})$ acts transitively on the elements of $\bar{V}$, then we will analyze $(\mathrm{Sp}(\bar{V}))_{y_{g}}$. Let $W$ be a symplectic submodule of $\bar{V}$ with basis $\{ x_1,\cdots,x_{g-1},y_1, \cdots, y_{g-1} \}$. Let also $J_0$ the induced symplectic form. That is,
	
	\[ J = \begin{pmatrix}
		J_0 & 0 \\
		0 & \begin{bmatrix} 0 & 1\\ -1 & 0 \end{bmatrix}
	\end{pmatrix} . \]
	
	Define $H(\bar{V}) = \{ (c,d) \mid c \in R^{2g-1} , \, d \in R \}$. We turn $H(\bar{V})$ into a group by the following product rule. For $(c,d), (c',d') \in H(\bar{V})$ we have
	
	\[ (c,d)(c',d') = ( c+c' , -c^T J_0 c'+d+d'). \]
	
	This is a well defined group since the product is a closed operation; the identity element is $(0,0)$ and the inverse of $(c,d)$ is $(-c,-d)$.
	
	\begin{proposition}
	There is a split short exact sequence
	
	\[ 1 \to H(\bar{V}) \to (\mathrm{Sp}(\bar{V}))_{y_{g}} \to \mathrm{Sp}(W) \to 1. \]
	\end{proposition}
	
	We will prove that every element of $(\mathrm{Sp}(\bar{V}))_{y_{g}}$ can be expressed by
	
	\[(A,c,d) \in \mathrm{Sp}(W) \times R^{2g-1} \times R.\]
	The right hand side of the expression above will be used in the proof of Theorem C.
	
	\begin{proof}
	Consider an element $S \in \mathrm{Sp}(\bar{V})$. We have
	
	\[ B = \begin{pmatrix}
	A & c & q \\
	p^T & a & e \\
	r^T & d & h
	\end{pmatrix}, \]
	where $c,q,p,r \in R^{2g-1}$, $a,e,d,h \in R$, and $A$ is a square matrix of dimension $g-1$. Since $By_g = y_g$, then $q=0$, $e=0$, and $h=1$. We have:
	
	\[ B^T J B = J \implies 
	\begin{pmatrix}
		A & p & r \\
		c^T & a & d \\
		0 & 0 & 1
	\end{pmatrix} \begin{pmatrix}
		J_0 & 0 \\
		0 & \begin{bmatrix} 0 & 1\\ -1 & 0 \end{bmatrix}
	\end{pmatrix} \begin{pmatrix}
		A & c & 0 \\
		p^T & a & 0 \\
		r^T & d & 1
	\end{pmatrix} = \begin{pmatrix}
	J_0 & 0 \\
	0 & \begin{bmatrix} 0 & 1\\ -1 & 0 \end{bmatrix}
	\end{pmatrix}
	\]
	
	By performing the calculations we end up to the system of equations
	
	\[
	\begin{matrix}
		A^T J_0 A - rp^t + pr^t &= J_0\\
		c^TJ_0 A - dp^T + ar^T &= 0\\
		-p^T &= 0\\
		A^T J_0 c - a r + dp &= 0\\
		c^T J_0 c &= 0\\
		a &= 1\\
		p &= 0
	\end{matrix}
	\]
	
	The equations above give $A^TJ_0A=J_0$, that is, $A \in \mathrm{Sp}(W)$. Also, $c^T J_0 c = 0$ and $c^TJ_0 A = -r^T$. The matrix $B$ has the following block form
	
	\[ B = \begin{pmatrix}
		A & c & 0 \\
		0 & 1 & 0 \\
		-c^TJ_0 A & d & 1
	\end{pmatrix}. \]
	
	Suppose that $B'$ is another matrix:
	\[ B' = \begin{pmatrix}
		A' & c' & 0 \\
		0 & 1 & 0 \\
		-c'^TJ_0 A' & d' & 1
	\end{pmatrix}. \]
	Then we have
	\[ B B' = \begin{pmatrix}
		AA' & Ac'+c & 0 \\
		0 & 1 & 0 \\
		-c^TJ_0 A A' -c'^TJ_0A' & -c^TJ_0Ac'+d+d' & 1
	\end{pmatrix}, \]
	where
	\[ -c^T J_0 A A' -c'^TJ_0A' = -c^T J_0 A A' -c'^T A^T J_0 A A' = -(c^T+c'^T A^T) J_0 A A'  \]
	and the latter one is equal to $-(Ac'+c)^T J_0 A A'$. Hence,
	\[ B B' = \begin{pmatrix}
		AA' & Ac'+c & 0 \\
		0 & 1 & 0 \\
		-(Ac'+c)^T J_0 A A' & -c^TJ_0Ac'+d+d' & 1
	\end{pmatrix}. \]
	Now the proposition follows easily by setting $B \mapsto (A,c,d)$.
	\end{proof}
	
	We note that since the rank of $\bar{V}$ is $2g$ over $R$, we write $\Sp{2g}{R}$ for $\mathrm{Sp}( \bar{V} )$ and we also write $(\Sp{2g}{R})_{y_g}$ instead of $(\mathrm{Sp}( \bar{V} ))_{y_g}$.

	\subsection{Basics on mapping class groups}
	
\begin{figure}[h]
		\begin{center}
			\includegraphics[scale=.4]{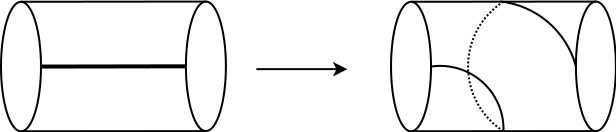}
		\end{center}
		\caption{The action of a Dehn twist $T_c$.}
		\begin{picture}(22,12)
			\put(-70,70){$\gamma$}
			\put(85,60){$T_c(\gamma)$}
			\put(0,70){$T_c$}
		\end{picture}
		\label{dt}
	\end{figure}

	Let $\Sigma_{g,n}$ be a surface of genus $g$ with $n$ boundary components. If $n=0$ then $\Sigma_g := \Sigma_{g,n}$. The \emph{mapping class group} $\Mod(\Sigma_{g,n})$ of $\Sigma_{g,n}$ is the group of isotopy classes of homeomorphisms that preserve the orientation of $\Sigma_{g,n}$ and fix the boundary pointwise. If $c$ is a simple closed curve of $\Sigma_{g,n}$, we denote by $T_c$ a Dehn twist along $c$. In this paper we make the convention that Dehn twists are right twists (see Figure \ref{dt}). We note that if two curves are homotopic then, their associated Dehn twists represent the same element in $\Mod(\Sigma_{g,n})$. \\
	
\begin{figure}[h]
		\begin{center}
			\includegraphics[scale=.5]{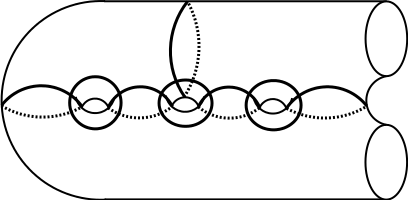}
		\end{center}
		\caption{The mapping class group is generated by Dehn twists about these curves.}
		\begin{picture}(22,12)
			\put(-22,78){$c_3$}
			\put(-56,73){$c_2$}
			\put(-78,75){$c_1$}
			\put(00,75){$c_4$}
			\put(22,75){$c_5$}
			\put(45,73){$c_6$}
			\put(73,76){$c_7$}
			\put(15,115){$b$}
		\end{picture}
		\label{mcg_gen}
	\end{figure}

	By work of Dehn the group $\Mod(\Sigma_{g})$ is finitely generated by Dehn twists \cite{dehn}. In Figure \ref{mcg_gen} we indicate the curves whose associated Dehn twists generate $\Mod(\Sigma_{g,n})$. By Humphries we have that $\Mod(\Sigma_{g})$ is generated by $T_{c_i}$, $T_b$ where $1 \leq i \leq 2g$ \cite{HU}. The latter result was extended by Johnson who proved that $\Mod(\Sigma_{g,1})$ is also generated by $T_{c_i}$, $T_b$ where $1 \leq i \leq 2g$ \cite{J1}. If $n=2$, then $\Mod(\Sigma_{g, n})$ is generated by $T_{c_i}$, $T_b$ where $1 \leq i \leq 2g+1$ \cite[Section 4.4.4]{BFM}.\\
	
	\subsection{Definition of the Torelli group}
	
	Here we define the Torelli group as a subgroup of $\Mod(\Sigma_{g,n})$ where $g > 0$ and $n\leq 2$. The definitions are generalized for any $g$ and $n$, but in this paper as explained in the introduction we do not require such a generality.
	
	\subsection*{Case $n \leq 1$.} For $n \leq 1$, then $\mathrm{H}_1(\Sigma_{g,n}) \cong \Zz^{2g}$ where generators of $\mathrm{H}_1(\Sigma_{g,n})$ are realized by simple closed curves in $\Sigma_{g,n}$. Let $\iota$ denote the algebraic intersection number between two curves in $\Sigma_{g,n}$. Then, $\iota$ is a symplectic form for $\mathrm{H}_1(\Sigma_{g,n})$ and we have a basis $\{ x_i, y_i \}$ as depicted in Figure \ref{sympb1} on the left. The action of $\Mod(\Sigma_{g,n})$ on $\mathrm{H}_1(\Sigma_{g,n})$ preserves the algebraic intersection number and hence we get a short exact sequence
	\[ 1 \to \I(\Sigma_{g,n}) \to \Mod(\Sigma_{g,n}) \to \mathrm{Sp}_{2g}(\Zz) \to 1. \]
	The group $\I(\Sigma_{g,n})$ is the Torelli group.\\
	
	Dehn twists in $\Mod(\Sigma_{g,n})$ are mapped to transvections in $\mathrm{Sp}_{2g}(\Zz)$. Since $\Mod(\Sigma_{g,n})$ is generated by finitely many Dehn twists, then $\mathrm{Sp}_{2g}(\Zz)$ is generated by finitely many transvections. We describe two basic types of elements in $\I(\Sigma_{g,n})$. Suppose that $\gamma \in \Sigma_{g,n}$ is a separating simple closed curve that bounds a subsurface of genus $h< g$ with one boundary component. The homology class of $\gamma$ is zero. Hence, the corresponding transvection in $\mathrm{Sp}_{2g}(\Zz)$ is a trivial transvection. Thus, $T_\gamma \in \I(\Sigma_{g,n})$ and it is called a \emph{Dehn twist about a separating simple closed curve of genus $h$}. Also, suppose that we have two simple closed curves $a,b$ in $\Sigma_{g,n}$ that bound a subsurface of $\Sigma_{g,n}$ with two boundary components and genus $h<g$. The corresponding homology classes of $a,b$ are equal. Hence, $T_a T_b^{-1}\in \I(\Sigma_{g,n})$. The latter mapping class is called \emph{bounding pair map of genus $h$}. Powell proved that these two types of elements generate $\I(\Sigma_{g,n})$ \cite{JP}. The result was improved by Johnson who proved that finite bounding pair maps suffice to generate $\I(\Sigma_{g,n})$ \cite{J1}.
	
	\subsection*{Case $n = 2$.} When $n \geq 2$ then $\mathrm{H}_1(\Sigma_{g,n}) \cong \Zz^{2g+n-1}$ is not longer symplectic. By computing the radical we get $\mathrm{rad}(\mathrm{H}_1(\Sigma_{g,n})) \cong \Zz^{n-1}$ with respect to the algebraic intersection number (alternating bilinear form). Furthermore, if we defined the Torelli group as the subgroup of $\Mod(\Sigma_{g,n})$ acting trivially on $\mathrm{H}_1(\Sigma_{g,n})$ we would face another problem. More precisely, an embedding $S \subset S'$ always induces a homomorphism $\Mod(S) \to \Mod(S')$. If we restrict to the Torelli groups on surfaces with many boundary components such homomorphisms do not necessarily hold. The definition we will give is due to Putman and he solves both the symplectic and the inclusion problem \cite{PM1}.\\
	
	We assume that $n=2$ but the definitions can be generalized for any $n$ \cite{PM1}. Let $S$ be a surface of genus $g$ with 2 boundary components, namely $\partial_1, \partial_2$. We choose a point $q_i$ in each $\partial_i$ and denote the set of such points by $Q$. By the long exact sequence of relative homology we get $\mathrm{H}_1(S,Q) \cong \Zz^{ 2g+2}$. \\
	
\begin{figure}[h]
		\begin{center}
			\includegraphics[scale=.45]{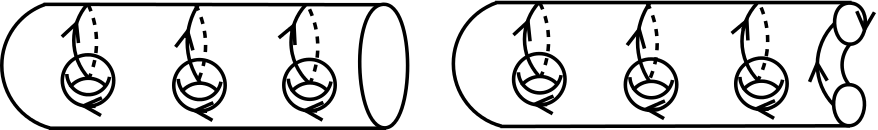}
		\end{center}
		\caption{Standard generators for $\mathrm{H}_1(\Sigma_{g,n})$, and $ \mathrm{H}_1(S,Q)$.}
		\begin{picture}(22,12)
			\put(-135,75){$y_1$}
			\put(-135,55){$x_1$}
			\put(-95,75){$y_2$}
			\put(-95,55){$x_2$}
			\put(-56,75){$y_3$}
			\put(-56,55){$x_3$}
			\put(30,75){$y_1$}
			\put(28,55){$x_1$}
			\put(71,75){$y_2$}
			\put(68,57){$x_2$}
			\put(108,75){$y_3$}
			\put(108,57){$x_3$}
			\put(174,78){$y_4$}
			\put(139,70){$x_4$}
		\end{picture}
		\label{sympb1}
	\end{figure}

	The module $\mathrm{H}_1(S,Q)$ is symplectic in the following manner. Consider the embedding $S \subset \Sigma_{g+1, 1}$ induced by gluing a pair of pants to the boundary components $\partial_1$ and $\partial_2$. As we saw earlier, $\mathrm{H}_1(\Sigma_{g+1,1};\Zz)$ is a symplectic module and actually, as a free module, it has rank $2g+2$. The inclusion $\mathrm{H}_1(S,Q) \to \mathrm{H}_1(\Sigma_{g+1,1};\Zz)$ is well defined and independent of the embedding \cite[Lemma 3.2]{PM1}. The symplectic form on $\mathrm{H}_1(\Sigma_{g+1,1};\Zz)$ restricts to a symplectic form on $\mathrm{H}_1(S,Q)$ and the basis of $\mathrm{H}_1(S,Q)$ is $x_i,y_i$ as indicated on the right hand side of Figure \ref{sympb1}.\\
	
	Note that $[\partial_1] = -[\partial_2]  \neq 0$ in $\mathrm{H}_1(S,Q)$ is a generator and generates an isotropic subspace of $\mathrm{H}_1(S,Q)$. Denote the latter module by $D$. Therefore, we have that
	
	\[ \mathrm{Aut}(\mathrm{H}_1(S,Q),D) \cong (\mathrm{Sp}_{2g+2}(\Zz))_{y_{g+1}}, \]
	
	where $y_{g+1}$ is the homology class of $\partial_1$ and 
	\[(\mathrm{Sp}_{2g+2}(\Zz))_{y_{g+1}} = \{ f \in \mathrm{Sp}_{2g+2}(\Zz) \mid f(y_{g+1})=y_{g+1} \}.\]
	We have that $\Mod(S)$ acts on $\mathrm{H}_1(S,Q)$ and this action induces a short exact sequence:
	
	\[ 1 \to \I(\Sigma_{g+1,1},S) \to \Mod(S) \to (\mathrm{Sp}_{2g+2}(\Zz))_{y_{g+1}} \to 1, \]
	
	where $\I(\Sigma_{g+1,1},S)$ is the Torelli group of $S$ such that $\Sigma_{g+1,1} \setminus S$ is a pair of pants. The surjection of the short exact sequence is due to Church \cite{TC} (see an earlier version in Arxiv arXiv:1108.4511v1). However, we can give an quick intuitive idea of the surjection.\\
	
	Consider the embedding $S \subset \Sigma_{g+1, 1}$ induced by the gluing a pair of pants to $\partial_1$ and $\partial_2$. Then $\partial_1$ is mapped to a nontrivial simple closed curve in $\Sigma_{g+1, 1}$ homologous to $y_{g+1}$ (see for example figure \ref{sympb1}). In terms of mapping class groups we have that the image of $\Mod(S) \to \Mod(\Sigma_{g+1,1})$ is the stabilizer of a simple closed curve representing $y_{g+1}$. The latter group acts on $\mathrm{H}_1(\Sigma_{g+1,1})$ and stabilizes $y_{g+1}$.\\
		
	If $\gamma$ is a simple closed curve in $S$ whose homology class in $\mathrm{H}_1(S,Q)$ is $[\gamma] = 0$, then $T_{\gamma}$ is called a Dehn twist about Q-separating curve. Also, if $a,b$ are two simple closed curves in $S$ such that $[a]=[b]$ in $\mathrm{H}_1(S,Q)$, then $T_aT_b^{-1}$ is called Q-bounding pair map. 
	
	Note that if $S \subset \Sigma_{g+1,1}$ is a surface with two boundary components where $\Sigma_{g+1,1} \setminus S$ is connected and $\gamma$ is a simple closed curve in $S\subset \Sigma_{g+1,1}$, then if the homology class of $\gamma$ is trivial in $\mathrm{H}_1(\Sigma_{g+1,1})$, this does not mean that the homology class of $\gamma$ is trivial in $\mathrm{H}_1(S,Q)$ too. We have the following theorem by Putman \cite[Theorem 1.3]{PM1}:
	
	\begin{theorem}[Putman]
		Let $S \subset \Sigma_{g+1.1}$ be a surface of genus $g$ with 2 boundary components, such that $\Sigma_{g+1.1} \setminus S$ is a pair of pants. The group $\I(\Sigma_{g+1,1},S)$ is generated by Dehn twists about Q-separating curves and Q-bounding pair maps.
	\end{theorem}
	
	Below we improve the result above in special cases.
	
	\begin{lemma}
		Let $S \subset \Sigma_{g+1.1}$ be a surface of genus $g\geq 3$ with 2 boundary components, such that $\Sigma_{g+1.1} \setminus S$ is a pair of pants. The group $\I(\Sigma_{g+1,1},S)$ is generated by Dehn twists about Q-bounding pair maps.
		\label{bp_maps}
	\end{lemma}
	
	The proof is similar to the Johnson's proof for $\I(\Sigma_{g,n})$ when $n \leq 1$ \cite{J2}.
	
	\begin{proof}
		The curves $x,y,z,b_1,b_2,b_3,b_4$ in Figure \ref{lant} form a lantern. Hence, we get
		\[ T_x T_y T_z = T_{b_1} T_{b_2} T_{b_3} T_{b_4} \implies T_{b_1} = (T_x T^{-1}_{b_4}) (T_y T^{-1}_{b_2}) (T_z T^{-1}_{b_3}) \]
		The elements $(T_x T^{-1}_{b_4}), (T_y T^{-1}_{b_2}), (T_z T^{-1}_{b_3})$ are Q-bounding pair maps and $T_{b_1}$ is a Dehn twist about a Q-separating curve. 
		
		\begin{figure}[h]
			\begin{center}
				\includegraphics[scale=.4]{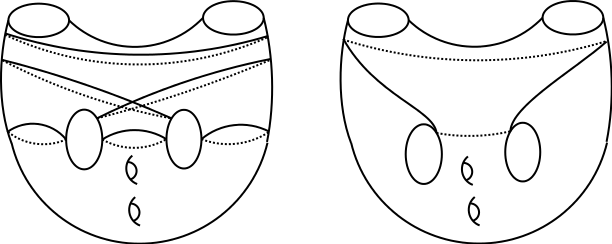}
			\end{center}
			\caption{Lantern formed by Q-bounding pairs and Q-separating curve.}
			\begin{picture}(22,12)
				\put(-120,110){$y$}
				
				\put(0,110){$x$}
				\put(0,123){$b_1$}
				
				\put(80,110){$z$}
				
				\put(-20,70){$b_2$}
				\put(-53,70){$b_3$}
				\put(-100,73){$b_4$}
			\end{picture}
			\label{lant}
		\end{figure}
		
		Note that $T_{b_1}$ is a separating curve of genus $g$. The result holds for all Dehn twists about separating simple closed curves of genus $g$ since they are all conjugate. For Dehn twists about separating simple closed curves of genus less than $g$, they are conjugate to a Dehn twist supported in a subsurface of $S$ bounded by a single separating curve. Since, this subsurface has one boundary component then, Dehn twists about separating curves are already products of bounding pair maps \cite{J2}.
	\end{proof}
	
	We can refine even more the generators of $\I(\Sigma_{g+1,1},S)$.
	
	\begin{lemma}
		The group $\I(\Sigma_{g+1,1},S)$ is generated by Q-bounding pair maps of genus 1.
		\label{bp_genus1}
	\end{lemma}
	
	\begin{proof}
		Bounding pair maps of genus less than $g$ are conjugate to bounding pair maps that lie in subsurfaces of $S$ bounded by a single separating curve. Hence, Lemma in these cases hold. There is a single bounding pair map of genus $g$, that is, $T_{\partial_1} T_{\partial_2}^{-1}$. We will see in the next section in Lemma \ref{johnson_relations} that $T_{\partial_1} T_{\partial_2}^{-1}$ is a product of bounding pair maps of genus less than $g$ (Relation (\ref{str3}) where $k=g+1$).
	\end{proof}
	
	\subsection{Abelianization of the Torelli group}
	
	The goal of this section is to describe the abelianization of the Torelli group $\I(\Sigma_{g+1,1},S)$. First we recall the abelianization of $\I(\Sigma_{g,1})$. In particular, we will describe two homomorphisms of $\I(\Sigma_{g,1})$ namely, the Johnson homomorphism $\tau$ \cite{J6} and the Birman-Craggs-Johnson homomorphism $\sigma$ \cite{J5} whose combination give the abelianization of $\I(\Sigma_{g,1})$ \cite{J3}. We will not construct the homomorphisms above, we will only give the formulas.\\
	
	\paragraph{One boundary component} Recall that $\mathrm{H}_1(\Sigma_{g,1})$ is generated by $x_i,y_i$, where $i \leq g$. Denote by $\Lambda^3_{g,1}$ the abelian group $\bigwedge^3 \mathrm{H}_1(\Sigma_{g,1})$ generated by $\alpha \wedge \beta \wedge \gamma$ where $\alpha, \beta, \gamma \in \{ x_i,y_i\}$. Since, $\mathrm{H}_1(\Sigma_{g,1}) \cong \Zz^{2g}$ we have that $\Lambda^3_{g,1}$ is torsion free of rank $\binom{2g}{3}$.\\
	
	Let $T_a T_b^{-1}$ be a bounding pair map of genus 1 bounding a subsurface $S' \subset \Sigma_{g,1}$. Choose $z,w \in \mathrm{H}_1(S') < \mathrm{H}_1(\Sigma_{g,1})$ such that $\iota(z,w)=1$ and let $c$ be the homology class of $a$ and $b$. Then,
	\[ \tau : \I(\Sigma_{g,1}) \to \Lambda_{g,1}^3 \]
	determined by $\tau(T_a T_b^{-1}) = z \wedge w \wedge c$. The action of $\Mod(\Sigma_{g,1})$ on $\I(\Sigma_{g,1})$ by conjugation induces an action of $\mathrm{Sp}_{2g}(\Zz)$ on $\Lambda_{g,1}^3$ making the latter abelian group as $\mathrm{Sp}$-module.\\
	
	Denote by $B_{g,1}^k$ the Boolean polynomial algebra $B^k(\mathrm{H}_1(\Sigma_{g,1};\Zz/2))$ of rank $k$. That is, $B_{g,1}^k$ is a $\Zz/2$-algebra (with 1) generated by $\overline{a}$ for each element $a \in \mathrm{H}_1(\Sigma_{g,1};\Zz/2)$ subject to relations:
	\begin{itemize}
		\item[1.] $\overline{a+b} = \overline{a}+\overline{b} + \iota(a,b) \bmod(\Zz/2)$,
		\item[2.] $\overline{a}^2 = \bar{a}$.
	\end{itemize}
	
	As a consequence $B_{g,1}^k$ has rank $\binom{2g}{k}+\binom{2g}{k-1}+\cdots+\binom{2g}{2}+\binom{2g}{1}+\binom{2g}{0}$. Let $T_a T_b^{-1}$ be a bounding pair map as above. Then, 
	\[ \sigma : \I(\Sigma_{g,1}) \to B_{g,1}^3 \]
	determined by $\sigma(T_a T_b^{-1}) = \bar{z} \bar{w} (\bar{c}+\bar{1})$.\\
	
	Let $a: \Lambda_{g,1}^3 \to \Lambda_{g,1}^3(\Zz/2)$ be the $\bmod(\Zz/2)$ reduction map, and let $b: B_{g,1}^3 \to \Lambda_{g,1}^3(\Zz/2)$ be such that $b(B_{g,1}^2) = 0$ and $b(\bar{\alpha} \bar{\beta} \bar{\gamma}) = \alpha \wedge \beta \wedge \gamma \bmod(\Zz/2)$. Define
	\[ W = \{ (u,v) \in \Lambda_{g,1}^3 \oplus B_{g,1}^3 \mid a(u) = b(v) \}. \]

	\begin{theorem}[Johnson \cite{J3}]
		Denote by $\I_{g,1}^{ab}$ the abelianization of $\I(\Sigma_{g,1})$. Then, we have an isomorphism
		
		\[ (\tau, \sigma): \I_{g,1}^{ab} \to W \]
		
		determined by $(\tau,\sigma)(\tilde{f}) = (\tau(\tilde{f}), \sigma(\tilde{f}))$, where $f \in \I(\Sigma_{g,1})$ and $\tilde{f} \in \I_{g,1}^{ab}$ is the corresponding element in the abelianization.
	\end{theorem}
	
	As a consequence, $\I_{g,1}^{ab}$ has rank $\sum_{i=0}^3 \binom{2g}{i}$.
	
	\paragraph{Two boundary components} In her thesis, van den Berg \cite{Berg}, calculated the abelianization of $\I(\Sigma_{g+1,1},S)$ \cite{Berg}. Van den Berg actually did it for any number of boundary components. We will describe the abelianization when $S$ has two boundary components and we will use it in the proof of Theorem C.\\
	
	Since $S \subset \Sigma_{g+1,1}$, then by definition we have $\I(\Sigma_{g+1,1},S) < \I(\Sigma_{g+1,1})$ and hence we deduce:
	
	\[
	\begin{matrix}
		\overline{\tau}: \I(\Sigma_{g+1,1},S) &\to \I(\Sigma_{g+1,1}) &\to \Lambda_{g+1,1}^3\\
		
		\overline{\sigma} :\I(\Sigma_{g+1,1},S) &\to \I(\Sigma_{g+1,1}) &\to B_{g+1,1}^3.
	\end{matrix}
	\]
	
	Define the module 
	\[ V_\Zz = \Zz x_1 \oplus \cdots \oplus \Zz x_g \bigoplus \Zz y_1 \oplus \cdots \oplus \Zz y_g \oplus \Zz y_{g+1} \]
	and $V_{\Zz/2}$ the mod(2) reduction of $V_\Zz$. Note that $V_\Zz < \mathrm{H}_1(\Sigma_{g+1,1};\Zz)$ and $V_{\Zz/2}< \mathrm{H}_1(\Sigma_{g+1,1};\Zz/2)$. Van den Berg proved the following theorem.
	
	\begin{theorem}
		The image of $\overline{ \tau }$ is $\bigwedge^3 V_\Zz$ and the image of $\overline{\sigma}$ is $B^3(V_{\Zz/2})$. Furthermore, if $\I_{g,2}^{ab}$ is the abelianization of $\I(\Sigma_{g+1,1},S)$, then 
		
		\[ (\tilde{\tau},\tilde{\sigma}): \I_{g,2}^{ab} \to \{ (u,v) \in \bigwedge^3 V_{\Zz} \bigoplus B^3(V_{\Zz/2}) \mid a(u)=b(v) \} \]
		
		is an isomorphism.
	\end{theorem}
	
	As a consequence, since $V_{R}$ (where $R$ is either $\Zz$ or $\Zz/2$) has rank $2g+1$, then $\I_{g,2}^{ab}$ has rank $\sum_{i=0}^3 \binom{2g+1}{i}$.
	
	\section{Chains maps and relations}
	
	In this section we provide a different perspective of bounding pair maps as given by Johnson. Then we form relations between some of them. These relations are helpful in the proof of Theorem A.
	
	\subsection{Chains}
	
	An \emph{chain} in $\Sigma_{g,n}$ is an ordered collection $(a_1,a_2,...,a_m)$ of simple closed curves with the following properties:
	\begin{itemize}
		\item[1.] The curves $a_i,a_{i+1}$ intersect transversely in a single point such that the algebraic intersection number between $a_i, a_{i+1}$ is $+1$.
		\item[2.] If $|i-j|>1$, then $a_i \cap a_j= \emptyset$.
		\item[3.] The homology classes of of $\{a_i\}$ in $\mathrm{H}_1(\Sigma_{g,n})$ are linearly independent.
	\end{itemize}
	
\begin{figure}[h]
		\begin{center}
			\includegraphics[scale=.65]{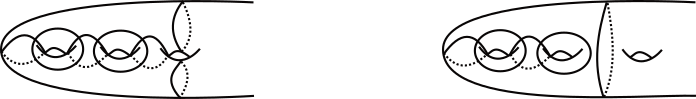}
		\end{center}
		\caption{Surfaces obtained by chains.}
		\begin{picture}(22,12)
			\put(-55,70){Positive curve $a$}
			\put(-55,50){Negative curve $a'$}

		\end{picture}
		\label{chains_surfaces}
	\end{figure}

	The length of an chain is the number of curves that it contains. We say that a chain is \emph{odd chain} if the number of curves in the chain is odd and \emph{even chain} otherwise. The length of a chain is the number of curves it contains. A regular neighborhood of an even chain of length $2k+2$ is (and $k\geq 1$) is a surface of genus $k$ with one boundary component, while a regular neighborhood of an odd chain of length $2k+1$ is a surface of genus $k$ with two boundary components. See for example Figure \ref{chains_surfaces}.

	Suppose we have a chain $(a_1, \cdots, a_{2k+1})$ with a regular neighborhood denoted by $M$. If we cut the $M$ along $a_1, a_3, \cdots, a_{2k+1}$ we obtain two subsurfaces, one lying on the left of $a_1, a_3, \cdots, a_{2k+1}$ (the positive subsurface) and the other lying on the right of $a_1, a_3, \cdots, a_{2k+1}$ (the negative subsurface) as in Figure \ref{chains_surfaces}. We call a positive curve $a$, the boundary of $M$ lying on the positive subsurface and we call a negative curve $a'$, the boundary of $M$ on the negative subsurface, as in Figure \ref{chains_surfaces}. 
	
	\begin{remark}
		This description above is ambiguous if we think any chain of curves on a surface, but the curves we are working on in this paper are $c_1$, $c_2$, ..., $b$ from Figure \ref{mcg_gen}.
	\end{remark}
	
	Consider an odd chain $(a_1,a_2,...,a_m)$, and let $a,a'$ be the curves of the boundary of a regular neighborhood of $(a_1,a_2,...,a_m)$ as above. Then we define
	\[ [a_1,a_2,...,a_m] = T_a T^{-1}_{a'}, \]
	and we will call it \emph{chain map}. In Figure \ref{chains_surfaces} on the left, we can see the chain $(c_1, c_2, c_3, c_4,c_5)$ and also by definition given above we have $T_aT_{a'}^{-1} = [c_1, c_2, c_3, c_4,c_5]$. If $f \in \Mod(\Sigma_{g,n})$, then $f*T_a T^{-1}_{a'}= f T_a T^{-1}_{a'}f^{-1} = T_{f(a)} T^{-1}_{f(a')}$. Likewise, we will write $f*[a_1,a_2,...,a_m]$ for $f[a_1,a_2,...,a_m]f^{-1} = [f(a_1),f(a_2),...,f(a_m)]$.\\
	
	\begin{figure}[h]
		\begin{center}
			\includegraphics[scale=.4]{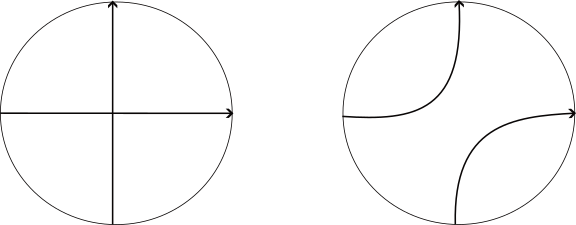}
		\end{center}
		\caption{Sum rule for curves of a single intersection point.}
		\label{chains_sum_rule}
	\end{figure}
	
	Let $a_i,a_{i+1}$ be two curves of an odd chain $(a_1,a_2,...,a_m)$. A regular neighborhood in the intersection point of $a_i,a_{i+1}$ looks like the left of Figure \ref{chains_sum_rule}. Then the new curve $a_i+a_{i+1} $ is produced by replacing the left side of \ref{chains_sum_rule} by the right side of \ref{chains_sum_rule} (see also Johnson's definition \cite[Section 2]{J1}. An \emph{odd subchain} of $(a_1,a_2,...,a_m)$ is a chain of the form $(k_1,k_2,...,k_l )$ such that $l$ is odd number with $l<m$, and $k_j = a_{i_j}+a_{i_{j}+1}+...+a_{i_{j+1}-1}, k_{j+1} =  a_{i_{j+1}}+a_{i_{j+1}+1}+...+a_{i_{j+2}-1}$.\\
	
	Consider the curves $c_i$ depicted in Figure \ref{mcg_gen}; consider also the odd chain $(c_1,c_2,...,c_{2g+1})$. Recall the curve $b$ as depicted in Figure \ref{mcg_gen}. Set $\beta = b+c_4$. Then, $T_b(c_4) = b+c_4 = \beta$.
	
	\begin{definition}
	Below we define two types of chains. A odd subchain of the form
	\[ (c_{i_1}+c_{{i_1}+1}+...+c_{{i_2}-1},c_{{i_2}}+...+c_{{i_3}-1},...,c_{{i_l-1}}+...+c_{{i_l}-1}) \]
	will be denoted by $(i_1, i_2,...,i_l)$, and the corresponding chain map by $[i_1, i_2,...,i_l]$. These chains will be called \emph{straight-chains}. An odd subchain of the form
	\[ (\beta + c_5 + \cdots + c_{{i_1}+1} + ... + c_{{i_2}-1}, c_{{i_2}} + ... + c_{{i_3}-1}, ..., c_{{i_l-1}} + ... + c_{{i_l}-1}) \]
	will be denoted by $(\beta, i_2,...,i_l)$, and the corresponding chain map by $[\beta, i_2,...,i_l]$. These chains will be called \emph{$\beta$-chains}.
	\label{strt_b_chains}
	\end{definition}
	
\begin{figure}[h]
		\begin{center}
			\includegraphics[scale=.65]{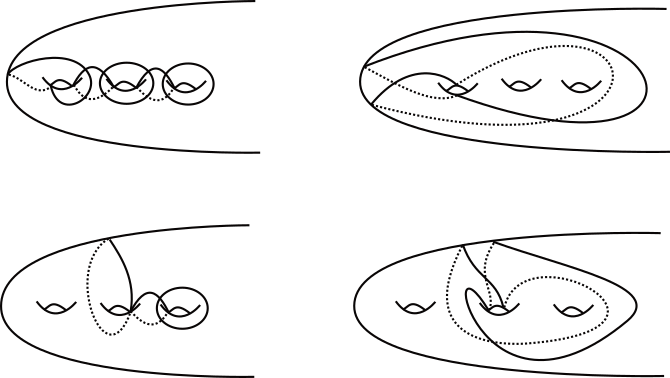}
		\end{center}
		\caption{On the top we have $[1,3,4,5,6] = T_a T_{a'}^{-1}$. On the bottom we have $[\beta, 5,6,7] = T_d T_{d'}^{-1}$}
		\begin{picture}(22,12)
			\put(-120,213){The chain $(1,3,4,5,6)$}
			\put(15,205){$a$}
			\put(15,180){$a'$}
			\put(-100,60){The chain $(\beta,5,6,7)$}
			\put(90,125){$d$}
			\put(70,125){$d'$}
		\end{picture}
		\label{examples_chains}
	\end{figure}

	For example $(c_1+c_2,c_3,c_4+c_5) = (1,3,4,6)$, $(c_1+c_2,c_3,c_4,c_5,c_6) = (1,3,4,5,6,7)$, $(\beta+c_5+c_6,c_7+c_8, c_9) = (\beta,7,9,10)$. In Figure \ref{examples_chains} we can see more examples of chains and chain maps. We have the following useful Lemma by Johnson \cite{J1}.
	
	\begin{lemma}
		The element $T_{c_j}$ commutes with $[i_1, i_2,...]$ if and only if $j,j+1$ are either both contained in or are disjoint from the $i's$. If $j=i_m$ but $j+1 \neq i_{m+1}$, then
		
		\begin{align}
			T_{c_{j}}^{-1}*[i_1, i_2,...] &= [i_1,...,i_{m-1},i_m+1,i_{m+1}...]. \label{A1} \tag{A1}\\
			T_{c_{j}}*[i_1 ,i_2,...] &= [i_1, i_2,...] [i_1,...,i_{m-1},i_m+1,i_{m+1},...]^{-1} [i_1, i_2...]. \label{A2} \tag{A2}
		\end{align}
		
		If $j+1=i_m$ but $j \neq i_{m-1}$, then
		\begin{align}
			T_{c_{j-1}}*[i_1, i_2,...] &= [i_1,...i_{m-1}, i_m-1, i_{m+1},...]. \label{B1} \tag{B1}\\
			T_{c_{j-1}}^{-1}*[i_1, i_2,...] &= [i_1, i_2,...] [i_1,...,i_{m-1}, i_m-1, i_{m+1},...]^{-1} [i_1, i_2,...]. \label{B2} \tag{B2}
		\end{align}
		\label{s_conj}
	\end{lemma}
	
	All statements of Lemma \ref{s_conj} except \ref{A2} and \ref{B2} are proven in \cite[Lemma 1]{J1} by Johnson. The proof of Equations \ref{A2} and \ref{B2} is given in the proof of \cite[Lemma 7]{J1}.
	
	\begin{definition}
		Let $G$ be a group and $H$ be a subgroup of $G$. We say that $g \in G$ \emph{normalizes} $h \in H$ in $H$ if $g h g^{-1} \in H$. Also, we say that $g\in G$ normalizes a subset $A \subset H$ in $H$ if for every $a\in A$, then $g a g^{-1} \in H$.
	\end{definition}
	
	Lemma \ref{s_conj} gives a consequence, given below.
	
	\begin{lemma}
		Consider a chain of curves $\{c_1, ...c_{2g}\}$ whose regular neighborhood is a surface of genus $g$ with one boundary component. Then, the group generated by the odd chain maps $[i_1, i_2, ..., i_{2g-1}]$ is normalized by $T^{\pm 1}_{c_i}$.
		\label{straight_norm}
	\end{lemma}
	
	The following theorem is due to Johnson \cite[Main Theorem]{J1}.
	
	\begin{theorem}
		Let $S$ be a surface of genus $g\geq 3$ with at most one boundary component. Then $\I(S)$ is generated by all odd straight chain maps and odd $\beta$-chain maps.
		\label{johns_ftheorem}
	\end{theorem}
	
	Before we proceed to the proof of Theorem A, we provide two relations proved by Johnson and they will be helpful in our proof.
	
	\begin{lemma}
	Let $k$ be a positive integer satisfying $3 \leq k \leq g$. We have the following relations.
	
	\begin{equation}
		\tag{J1}
		[2,3,4,...,2k+1] = f [4,5,6,...,2k+1] [2,3,4,5]
		\label{str1}
	\end{equation}
	
	where $f$ is a bounding pair map of genus 1.
		
	\begin{equation}
		\tag{J2}
		[2,3,4,...,2k+1]^{-1} (T_b)^*[2,3,4,...,2k+1] = [2,3,4,5]^{-1} [4,5,6,...,2k+1]^{-1} [\beta,5,6,...,2k+1] (T_b)^*[2,3,4,5]
		\label{str2}
	\end{equation}
	
	\begin{equation}
		\tag{J3}
		[1,2,3,4] [1,2,5,6,...,2k] (T_b)^*[3,4,5,...,2k] = [5,6,...,2k] [1,2,3,4,...,2k]
		\label{str3}
	\end{equation}
		\label{johnson_relations}
	\end{lemma}

	\subsection{Proof of Theorem A}
	
	Let $\Sigma_{g+1,1}$ be a surface of genus $g+1$ with one boundary component. Let $S \subset \Sigma_{g+1,1}$ be a subsurface of genus $g$ with two boundary components, say $\partial_1, \partial_2$, such that $\Sigma_{g+1, 1} \setminus S$ is a pair o pants (recall from Figure \ref{surface_model}). To prove Theorem A we need to find a finite generating set for $\I(\Sigma_{g+1,1},S)$. Recall the chain maps described in the previous subsection with examples in Figures \ref{chains_surfaces} and \ref{examples_chains}. Our aim is to prove the following Theorem.
	
	\begin{theorem}
		For $g\geq 3$ the group $\I(\Sigma_{g+1,1},S)$ (where the genus of $S$ is $g$) is generated by the odd sub-chains of $[i_1, i_2, ..., 2g+2]$ and $[\beta, i_1, i_2, ..., 2g+2]$.
		\label{first_th}
	\end{theorem}
	
	Recall that the embedding $S \subset \Sigma_{g+1,1}$ induces an injection $\iota: \I(\Sigma_{g+1,1},S) \to \I(\Sigma_{g+1,1})$. We can see that if $f \in \I(\Sigma_{g+1,1},S)$ is a generator from Theorem \ref{first_th}, then $\iota(f)$ is a generator of $\I(\Sigma_{g+1,1})$ provided in Theorem \ref{johns_ftheorem} (as described in the introduction and in Figure \ref{surface_model}).\\ 
	
	To prove Theorem \ref{first_th} we proceed as follows. First of all, recall that $S$ is homeomorphic to $\Sigma_{g,2}$ and by Section 2, $\I(\Sigma_{g+1,1},S)$ is normal in $\Mod(S)$. Suppose that $J_g < \I(\Sigma_{g+1,1},S)$ is generated by the bounding pair maps of Theorem \ref{first_th}. If every generator of $\Mod(S)$ normalizes every generator of $J_g$, then $J_g$ is normal in $\Mod(S)$. But since $\I(\Sigma_{g+1,1},S)$ is normally generated by Q-bounding pair maps, then by Lemma \ref{bp_maps} we deduce that $J_g = \I(\Sigma_{g+1,1},S)$.\\
	
	Before we proceed to the proof of Theorem \ref{first_th} we need the following useful lemmas.
	
	\begin{lemma}[Normalization trick]
		Let $G$ be a group and $H$ be a subgroup of $G$. Consider the elements $f,g \in G$ and for every $p \in H$, we have $g*p \in H$ and $[f,g] = 1.$ Then, $f*p=f p f^{-1} \in H$ if and only if $fg*p=f*(g*p) \in H$.
		\label{norm_trick}
	\end{lemma}
	
	\begin{proof}
		Suppose that $f (g p g^{-1}) f^{-1} \in H$. Then,
		\[ f (g p g^{-1}) f^{-1} = g (f p f^{-1}) g^{-1} = g g_1 g^{-1} \in H. \]
		Since $g g_1 g^{-1} \in H$, then the latter element is expressed as a product of elements of $H$. That is, $g g_1 g^{-1} = h_1 h_2 \cdots h_k$. Then, $g_1 = g^{-1} h_1 h_2 \cdots h_k g$ and the latter element is in $H$ since $g$ normalizes $H$. Hence, $g_1 \in H$, that is $f p f^{-1} \in H$.\\
		
		Suppose now that $f p f^{-1} \in H$. Then $g(f p f^{-1})g^{-1} \in H$ by hypothesis. But $g(f p f^{-1})g^{-1} = f(g p g^{-1})f^{-1}$.
	\end{proof}
	
	\begin{figure}[h]
		\begin{center}
			\includegraphics[scale=.6]{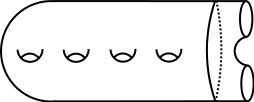}
		\end{center}
		\caption{Embedding of $\Sigma_{g,1}$ to $S$.}
		\label{oneboundarysubsurface}
	\end{figure}
	
	The proof of Theorem \ref{first_th} will be proved in four Lemmas. Recall the generators of $\Mod(S)$ as depicted in Figure \ref{mcg_gen}. In the first Lemma below we prove that $T^{\pm 1}_{c_j}$ normalize the generators of $J_g$ in $J_g$ if $j\leq 2$ or $j \geq 5$. Actually, this Lemma gives a stronger result (see below for details). In the next three Lemmas we prove that $T^{\pm 1}_{c_3}$, $T^{\pm 1}_{c_4}$, $T^{\pm 1}_b$ normalize the generators of $J_g$ in $J_g$.
	
	\begin{lemma}
		Let $T_{c_j}$ where $j \leq 2g+1$ be the generators of $\Mod(S)$ as indicated in Figure \ref{mcg_gen}.
		\begin{itemize}
			\item If $p$ is a straight chain map, then $T_{c_j}^{\pm 1}*p \in J_g$ for all $c_j$.
			\item If $p$ is a beta chain map, then $T_{c_j}^{\pm 1}*p \in J_g$ for $c_j$ when either $j \geq 5$ or $j \leq 2$.
		\end{itemize}
		\label{normal_argument}
	\end{lemma}
	
	\begin{proof}
		Suppose that $p=(i_1,\cdots,i_{2k})$ be a straight chain. By Lemma \ref{s_conj}, $T_{c_j}^{\pm 1}$ commutes with $p$ if $j,j+1$ are either both contained or a disjoint from the indices of $p$. If not, enlarge $p$ by including both $j,j+1$ and denote the new chain by $p'$. The regular neighborhood of $p'$ is homeomorphic to a surface with one boundary component $M$, as indicated in Figure \ref{chains_surfaces} on the right. Hence, by Lemma \ref{straight_norm} $T_{c_j}^{\pm 1}*p$ is a product of odd maximal straight chain maps. These chain maps are elements of $J_g$. Hence, $T_{c_j}^{\pm 1}*p \in J_g$.\\
		
		Suppose that $p=(\beta,i_2,\cdots, i_{2k})$. If $j\leq 2$, then $c_j$ is disjoint from the curves of $p$, hence $T_{c_j}^{-1}$ commutes with $p$.\\
		
		When $j\geq 5$ we apply the same argument we did for straight chains as follows. First of all, $T_{c_j}^{\pm 1}$ commutes with $p$ if $j,j+1$ are either both contained or a disjoint from the indices of $p$. If not, enlarge $p$ by including both $j,j+1$ and denote it by $p'$. The new chain $p'$ is a chain starting from $\beta$ and the regular neighborhood is a surface with one boundary component $M$ (again homeomorphic to the surface in Figure \ref{examples_chains} on the left). By Lemma \ref{straight_norm}, the Dehn twists $T_{\beta}^{\pm 1}$, $T_{c_5}^{\pm 1}, T_{c_6}^{\pm 1} \cdots $ normalize maximal odd chain maps of $M$. However, these chain maps are $\beta$-chain maps in $S$ and consequently in $J_g$. Since $p$ is one of these $\beta$-chain maps, then, $T_{c_j}^{\pm 1}*p \in J_g$ if $j\geq 5$.
	\end{proof}
	
	Consider the embedding $\Sigma_{g,1} \subset S$ as indicated in Figure \ref{oneboundarysubsurface}. The main strategy in the proof below is to conjugate an chain map by appropriate elements from Lemma \ref{normal_argument} of $\Mod(S)$ so that the result is inside $\I(\Sigma_{g,1})<J_g$.
		
	\begin{lemma}
		Suppose that $p$ is a generator of $J_g$. Then, $T^{\pm 1}_{c_3}*p \in J_g$.
		\label{normalize_C3}
	\end{lemma}
	
	\begin{proof}
		By the first bullet of Lemma \ref{normal_argument}, then if $p$ is a straight-chain of $J_g$, then $T^{\pm 1}_{c_3}*p \in J_g$. Suppose that $p$ is a $\beta$-chain. Definitely it is true that
		\[ T^{\pm 1}_{c_3}*[\beta, \cdots, 2g+1] \in \Mod(\Sigma_{g,1}) \subset J_g. \]
		Suppose that $p = [\beta, i_1, i_2, \cdots, 2g+2]$. The length of $p$ is at most $2g-3$ and it is impossible to find consecutive indices $[\beta,5,6,\cdots,2g+2]$. Because that would exceed the length $2g-3$. Hence, there exist at least two consecutive indices $i_k, i_{k+1}$ in $p$ such that $i_{k+1}-i_k\geq 2$. Therefore, there exists a product of Dehn twists $T_{c_j}^{\pm 1}$, where $j \geq i_k$ with the following property. Call this product by $f$. Then $[f, T_{c_3}] = 1$ and
		\[ f*[\beta, i_1, i_2, \cdots, i_k, i_{k+1} \cdots, 2g+2] = [\beta, i_1, i_2, \cdots, i_k, i_k+1, i_k+2, \cdots, 2g+1] \in \Mod(\Sigma_{g,1}) \subset J_g. \]
		By the normalization trick we deduce that $T^{\pm 1}_{c_3}*p \in J_g$.
	\end{proof}
	
	\begin{lemma}
		Suppose that $p$ is a generator of $J_g$. Then, $T^{\pm 1}_{c_4}*p \in J_g$.
		\label{normalize_C4}
	\end{lemma}
	
	\begin{proof}
	Using Lemma \ref{normal_argument}, then $T^{\pm 1}_{c_4}$ normalizes all straight chain maps in $J_g$. As in Lemma \ref{normalize_C3}, we want to prove $T^{\pm 1}_{c_4}*[\beta, i_1, i_2,\cdots, 2g+2]$. If $p = [\beta, 6, \cdots, 2g+2]$, then $T^{\pm 1}_{c_4}$ commutes with $p$ and we get $T^{\pm 1}_{c_4}*p=p$. If $p = [\beta, 5, \cdots, 2g+2]$, then like in the proof of Lemma \ref{normalize_C3} there exists $i_k,i_{k+1}$ with $i_{k+1}-i_k\geq 2$ and $i_k>5$. Hence, there is a product of $T^{\pm 1}_{c_j}$ with $j\geq i_k$, say $f$, such that if we conjugate $p$ by this $f$ we get
	\[ f*[\beta, 5, \cdots, i_k , i_{k+1}, \cdots 2g+2] = [\beta, 5, \cdots, i_k , i_k+1, i_k+2, \cdots 2g+1] \in J_g. \]
	
	By the normalization trick we deduce that $T^{\pm 1}_{c_4}*p \in J_g$.
	\end{proof}
	
	\begin{lemma}
		
		Suppose that $p$ is a generator of $J_g$. Then, $T^{\pm 1}_b*p \in J_g$.
		\label{normalize_b}
	
	\end{lemma}
	
	\begin{proof}
		We want to prove that $T^{\pm 1}_b$ normalize straight-chain maps and $\beta$-chain maps in $J_g$. Recall that every chain map $[\cdots, 2g+1] \in \I(\Sigma_{g,1})<J_g$ is normalized by $T^{\pm 1}_b$. We have the following chain maps that we want to conjugate by $T^{\pm 1}_b$.
		\begin{align}
		[\beta,\cdots, 2g+2] \label{CM1} \\
		[4,\cdots, 2g+2] \label{CM2} \\
		[i,4,\cdots,2g+2] \label{CM3} \\
		[i_1,i_2,4,\cdots, 2g+2] \label{CM4} \\
		[1,2,3,4,\cdots,2g+2] \label{CM5}
		\end{align}
		
		The chain maps of type \ref{CM5} are disjoint from $T_b^{\pm 1}$, hence, $T_b^{\pm 1}*[1,2,3,4,\cdots,2g+2]=[1,2,3,4,\cdots,2g+2]$.\\
				
		The chain maps of types \ref{CM1} and \ref{CM2} have maximal length $2g-3$. The chain maps of type \ref{CM4} have maximal length $2g-1$. This means that there is not a chain map of the form
		\[ [\cdots, 4,5,6,\cdots 2g+1,2g+2], \]
		that is, consecutive indices from 4 up to $2g+2$. Because the length of this kind of chain would exceed the maximal lengths indicated above. Hence, there exist $f$ which is a product of $T_{c_j}^{\pm 1}$ where $j\geq 5$ and
		\[ f*[\cdots, 4,\cdots, 2g+2] = [\cdots, 4,5,6,\cdots, 2g, 2g+1] \in J_g. \]
		By the normalization trick we are done for these cases. \\
		
		Finally, we want to prove that $T_b^{\pm 1}$ normalizes chain maps of type \ref{CM3}. Since $T_{c_1}^{\pm 1}, T_{c_2}^{\pm 1}$ normalize all generators of $J_g$ in $J_g$ by Lemma \ref{normal_argument} and using the normalization trick, we only need to prove that
		\[ T_b^{\pm 1}* [3,4,\cdots, 2g+2] \in J_g. \]
		Setting $k=g+1$ in Relation \ref{str3} we deduce that $T_b* [3,4,\cdots, 2g+2] \in J_g$. We conjugate the same relations by $f = T_{c_2}^{-1} T_{c_3}^{-1}T_{c_1}^{-1}T_{c_2}^{-1}T_b^{-1}$ to deduce
		\[ f*[1,2,3,4] T_b^{-1}*[3,4,5,\cdots 2g+2] f*[3,4,5,\cdots 2g+2] = [5,6,\cdots 2g+2] [1,2,3,4,\cdots , 2g+2]. \]
		We have that $f*[1,2,3,4] \in J_g$ by \cite[Lemma 6]{J1}. Also, $f*[3,4,5,\cdots 2g+2] \in J_g$ by the first bullet of Lemma \ref{normal_argument}. Therefore, $T_b^{- 1}* [3,4,\cdots, 2g+2] \in J_g$.
	\end{proof}
	
	\begin{corollary}
		If the genus of $S$ is 3, then  $\I(\Sigma_{g+1, 1},S)$ is generated by 103 elements.
		\label{numerical_thmA}
	\end{corollary}
	
	\begin{proof}
		If the genus of $S$ is 3, then  $\I(\Sigma_{g+1, 1},S)$ is generated by
		\begin{align*}
			[1,2,3,4,5,6,8]\\
			[i_1,i_2,i_3,i_4]\\
			[i_1,i_2,i_3,i_4,i_5,i_6]\\
			[\beta,j_1,j_2,j_3]
		\end{align*}
		where $1 \leq i_1 < i_2 < \cdots \leq 8$ and $1 \leq j_1 < i_2 < j_3 \leq 8$. Therefore, $\I(\Sigma_{g+1, 1},S)$ is generated by
		\[ 1 + \binom{8}{4} + \binom{8}{6} + \binom{4}{3} = 103 \]
		elements.
	\end{proof}
	
	\section{Small generating sets for the Torelli group}
	\label{small_gen_sets}
	
	In this Section we provide different generating sets for the Torelli group, other than Theorem \ref{first_th}. This will prove Theorem B. Let $S \subset \Sigma_{g+1,1}$ be a surface of genus $g$ with two boundary components, such that $\Sigma_{ g+1,1} \setminus S$ is connected. Put $\boldsymbol{g} = \{1,2,\cdots g\}$ and let $I \subseteq \boldsymbol{g}$. Below we will define subsurfaces $S_I\subset S \subset \Sigma_{g+1,1}$ for any $I \subseteq \boldsymbol{g}$ such that 
	\begin{itemize}
		\item $|I| \leq g$,
		\item the genus of $S_I$ is $|I|$,
		\item $S \setminus S_I$ is connected.
		\item $S$ has two boundary components.
	\end{itemize}
	
	The main Theorem of this section is to prove that $\mathcal{I}(\Sigma_{g+1,1},S)$ is generated by the union of all $\mathcal{I}( \Sigma_{g+1,1},S_I)$. Theorem B is then implied.
	
	\subsection{Subsurfaces of lower genus}
	
	\label{s_i}
	
	Let $S$ be a surface of genus $g$ with two boundary components and suppose that $S\subset \Sigma_{g+1,1}$ such that $\Sigma_{g+1,1} \setminus S$ is connected.\\
	
	\begin{figure}[h]
		\begin{center}
			\includegraphics[scale=.25]{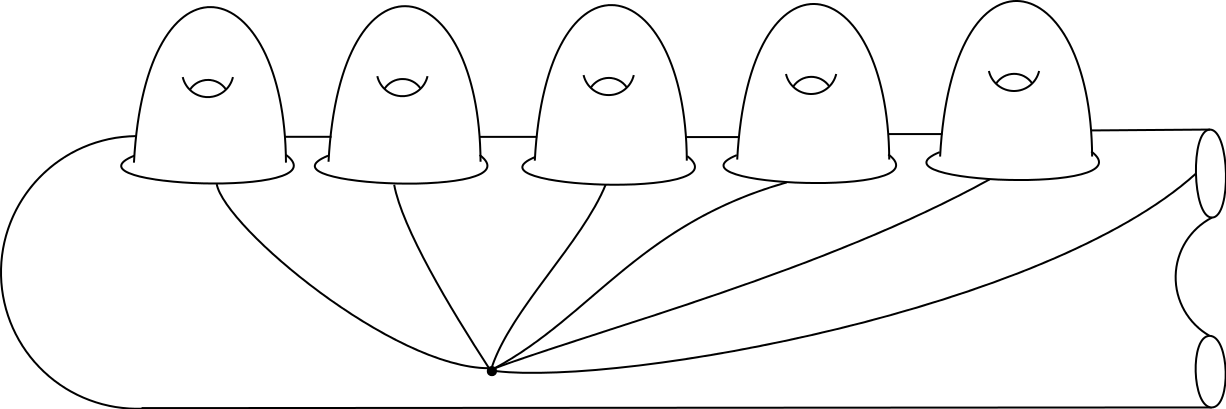}
			\includegraphics[scale=.25]{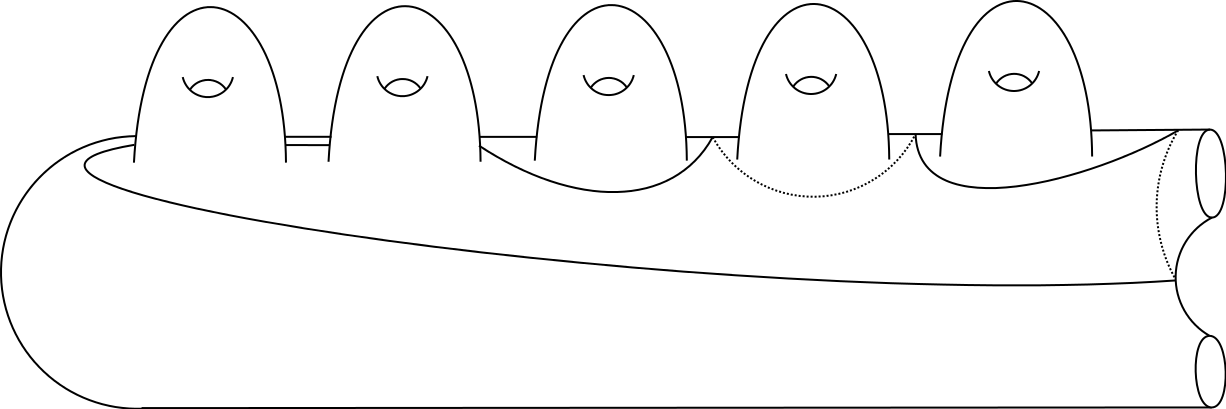}
		\end{center}
		\caption{A surface with two boundary components with an example of a subsurface $S_I$.}
		\begin{picture}(22,12)
			\put(-90,190){$\delta_1$}
			\put(-95,215){$R_1$}
			\put(-50,190){$\delta_2$}
			\put(-50,215){$R_2$}
			\put(-15,193){$\delta_3$}
			\put(0,215){$R_3$}
			\put(36,195){$\delta_4$}
			\put(55,215){$R_4$}
			\put(80,190){ $\delta_5$}
			\put(100,215){ $R_5$}
			\put(90,170){ $\delta_6$}
			\put(170,210){$\partial_1$}
			\put(170,160){$\partial_2$}
			\put(-30,160){$p_0$}
			\put(36,95){$S_{ \{1,2,4\} }$}
		\end{picture}
		\label{partitioned}
	\end{figure}
	
	Now we will define subsurfaces $S_I$ for any $I \subset \boldsymbol{g }$. For each $i\leq g$ denote by $R_i$ the subsurfaces of $S$ homeomorphic to a torus with one boundary component as depicted in Figure \ref{partitioned}. Choose a point $p_0$ in the interior of $S'$. Fix a point $*_i$ in each $\partial R_i$ for $i \leq g$ and $*_{g+1}$ in $\partial_1$. For $i\geq 1$ let $\delta_i$ be arcs such that
	\begin{itemize}
		\item $\partial \delta_i = \{ p_0, *_i \}$ with orientation that starts from $p_0$ and ends to $*_i$.
		\item $\delta_i \cap \partial R_j = \emptyset$ if $i \neq j$.
		\item $\delta_i \cap \partial_1 = \emptyset$ if $i \leq g$.
		\item $\delta_i \cap \partial  R_i = *_i$.
		\item $\delta_i \cap \partial_1 = *_i$ if $i = g+1$.
		\item $\delta_i \cap \delta_j = p_0$ if $i\neq j$.
	\end{itemize}
	
	Consider the subset $I = \{ i_1, \cdots, i_k \} \subseteq \boldsymbol{g}$ and denote by $S_I$ the regular neighborhood of $p_0$, $\delta_{ i_1}$, $R_{i_1}$, $\cdots$, $\delta_{ i_k}$, $R_{i_k}$, $\delta_{ g+1}$, $\partial_1$. We have that $|I| = k$ and the genus of $S_I$ is $k$ with two boundary components, where one of them is $\partial_1$. In Figure \ref{partitioned} we show an example where $I = \{ 1,2,4 \}$. Since $S_I \subset S \subseteq \Sigma_{g+1,1}$ and $\Sigma_{g+1,1} \setminus S_I$ is connected, then we can define $\mathcal{I}(\Sigma_{g+1,1},S_I)$. Hence $\mathcal{I}(\Sigma_{g+1,1},S_I) < \mathcal{I}(\Sigma_{g+1,1},S)$.
	
	\begin{remark}
		The subsurfaces $S_I$ are related to the subsurfaces of \cite[Section 4.2]{PM3} by Church-Putman as follows. If we cap the boundary $\partial_1$ in $S$, then our definitions above imply the definitions of \cite[Section 4.2]{PM3}.
	\end{remark}
	
	\subsection{Calculation of smaller generating sets}
	
	Theorem B is a consequence of Theorem \ref{cubic} below. For $I \subset \boldsymbol{g}$, let $S_I$ be the subsurfaces defined above. The next Theorem proves Theorem B.
	
	\begin{theorem}
		Let $S \subset \Sigma_{g+1,1}$ be a genus $g\geq 3$ surface with two boundary components such that $\Sigma_{g+1,1} \setminus S$ is connected and define $S_I \subset S$ as in Section \ref{s_i}. Fix a positive integer $k\geq 3$. Then
		\[ \mathcal{I}(\Sigma_{g+1,1},S) = \langle \mathcal{I}(\Sigma_{g+1,1},S_I) \mid I \subset \boldsymbol{g} \, \text{satisfies} \, |I| = k \rangle. \]
		\label{cubic}
	\end{theorem}
	
	Theorem \ref{cubic} is the analogue of a result by Church-Putman \cite[Proposition 4.5]{PM3}. Also, the proof of Theorem \ref{cubic} is similar as in the proof of \cite[Proposition 4.5]{PM3}. More precisely, \cite[Proposition 4.5]{PM3} is based on \cite[Theorem B]{PM2}. The next lemma is the analogue of \cite[Theorem B]{PM2}.
	
	\begin{lemma}
		Let $\Sigma_{g+1,1}, S$ be as above. For $I \subset \boldsymbol{g}$  denote by $X_I$ the closure of
		\[ S \setminus \bigcup_{l \notin I   } R_l. \]
		Then, $\mathcal{I}(\Sigma_{g+1,1},S)$ is generated by
		\[ \bigcup_{|I|=3} \mathcal{I}(\Sigma_{g+1,1}, X_I). \]
		\label{cut_handles}
	\end{lemma}
	
	The subsurface $X_I \subset S \subset \Sigma_{g+1,1}$ is homeomorphic to a surface of genus $|I|$ with $2+g-|I|$ boundary components. If we cap the boundaries $\partial_1,\partial_2$ and we restrict to $|I|=3$, then we obtain the generating set of \cite[Theorem B]{PM2}. To prove Lemma \ref{cut_handles} we need first the following Lemma (following the strategy in \cite{PM2}).
	
	\begin{figure}[h]
		\begin{center}
			\includegraphics[scale=.6]{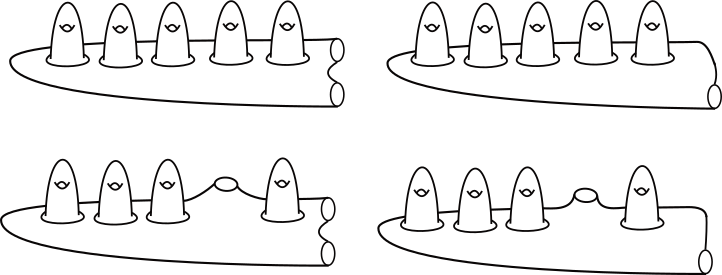}
		\end{center}
		\caption{On the top left we have the surface $S$. On the top right we have $S'$. On the bottom left we have $Y_4$ and on the bottom right we have $Y^1_4$. We can see that $Y_4 \subset S \subset S'$ and $Y_4 \subset Y^1_4 \subset S'$.} 
		\begin{picture}(22,12)
			\put(-50,145){$S$}
			\put(-50,80){$Y_4$}
			\put(115,75){$Y_4^1$}
			\put(115,145){$S'$}
		\end{picture}
		\label{y_subsurfaces}
	\end{figure}
	
	\begin{lemma}
		Let $\Sigma_{g+1,1}, S$ be as above. Define $Y_i$ as the closure of $S \setminus R_i$. The group $\mathcal{I}( \Sigma_{g+1,1},Y_i)$ is generated by
		\[ \bigcup_{1 \leq i \leq g} \mathcal{I}(\Sigma_{g+1,1}, Y_i). \]
		\label{handle_torelli}
	\end{lemma}
	
	We note that $Y_i \subset S$ is homeomorphic to a surface of genus $g-1$ with 3 boundary components, where two of them are $\partial_1, \partial_2$ and the third is $\partial R_i$.
	
	\begin{proof}[Proof of Lemma \ref{handle_torelli}]
		
		We divide the proof into two steps.
		
		\paragraph{Step 1} Let $S'$ be a surface of genus $g\geq3$ with 1 boundary component obtained from $S$ after capping $\partial_1$. Note that $S \Subset S'$ and $S'$ is not a subsurface of $\Sigma_{ g+1,1}$.
		
		Denote by $Y_i^1$ the surface obtained by $Y_i$ after capping $\partial_1$. Observe that $Y_i^1$ is a subsurface of $S'$. The surface $Y_i^1$ is homeomorphic to a surface of genus $g-1$ with 2 boundary components (one is $\partial_2$ and the other is $\partial R_i$). We will show that $\I(S')$ is generated by
		\[ \bigcup_{1 \leq i \leq g} \I(S' , Y_i^1). \]
		
		Consider the Birman exact sequence:
		\[ 1 \to  \pi_1(U \overline{S}) \to \I(S') \to \I(\overline{S}) \to 1, \]
		where $\overline{S}$ is a closed surfaces of genus $g-1$ obtained from $S'$ after capping $\partial_2$ and $U\Sigma_{g}$ is the unit tangent bundle of $\overline{S}$. We consider $\pi_1(U \overline{S})$ as the disk-pushing map subgroup (see Farb-Margalit's book \cite{BFM} for more details) with the following generators.
		
		\begin{figure}[h]
			\begin{center}
				\includegraphics[scale=.6]{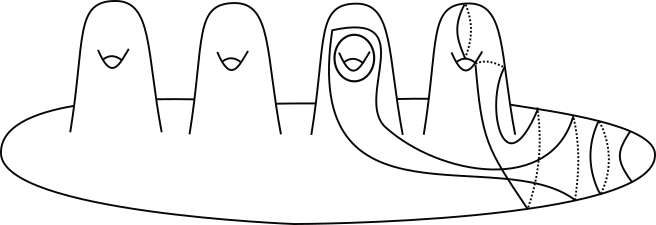}
			\end{center}
			\caption{The generators of the disk pushing map in $\I(S')$.} 
			\begin{picture}(22,12)
				\put(20,96){$\tilde{x}_3$}
				\put(10,70){$\tilde{x}'_3$}
				\put(73,150){$\tilde{y}_4$}
				\put(70,95){$\tilde{y}'_4$}
			\end{picture}
			\label{disk_pushing_map}
		\end{figure}
		
		Denote by $\tilde{x}_i, \tilde{y}_i$ the longitude and the meridian of $R_i$ as depicted in Figure \ref{disk_pushing_map}. Construct a simple closed curve $\tilde{x}'_i$ as follows. Glue a disk $d$ in the boundary $\partial_2$. Then drag a segment of $\tilde{x}_i$ along the surface until it crosses the disk $d$. Remove the disk $d$. The result is $\tilde{x}'_i$. Likewise, we construct $\tilde{y}'_i$. There are examples of $\tilde{x}_i$, $\tilde{x}'_i$, $\tilde{y}_i$, $\tilde{y}'_i$ in Figure \ref{disk_pushing_map}. Denote also by $D$ the separating simple closed curve parallel to $\partial_2$. The group $\pi_1(U \overline{S} )$ is generated by $T_{\tilde{x}_i} T_{\tilde{x}_i}^{-1}$, $T_{\tilde{y}_i} T_{\tilde{y}_i}^{-1}$, $T_D$. We can see that each generator of $\pi_1(U \overline{S} )$ lies in $\bigcup \mathcal{I}(S', Y_i^1)$. Also, via the Birman exact sequence the set $\bigcup_{1 \leq i \leq g} \mathcal{I}(S', Y_i^1)$ is mapped to a generating set of $\I(\overline{S})$ (Denote by $\tilde{Y}_i^1$ the surface obtained by $Y_i^1$ after capping $\partial_2$. Putman in the proof of \cite[Theorem B]{PM2} proved that $\I(\overline{S})$ is generated by $\bigcup_{1 \leq i \leq g} \I(\overline{S},\tilde{Y}_i^1)$. This is exactly the image of $\bigcup_{1 \leq i \leq g} \mathcal{I}(S', Y_i^1)$ in $\I(\overline{S})$). Hence, $\mathcal{I}(S')$ is generated by 
		\[\bigcup_{1 \leq i \leq g} \mathcal{I}(\Sigma_{g+1,1}, Y_i^1).\]
		
		\paragraph{Step 2} Let $S$ be a surface of genus $g\geq 3$ with 2 boundary components. Consider the Birman exact sequence:
		
		\[ 1 \to [\pi_1(\tilde{S}), \pi_1(\tilde{S})] \to \mathcal{I}(\Sigma_{g+1,1},S ) \to \mathcal{I}(\tilde{S} ) \to 1, \]
		
		where $\tilde{S}$ is obtained by $S$ after capping a boundary component. We have that $\bigcup_{1 \leq i \leq g} \mathcal{I}(\Sigma_{g+1,1}, Y_i)$ is mapped to the generating set of $\mathcal{I}(\tilde{S} )$ by Step 1. Also, $[\pi_1(\tilde{S}), \pi_1( \tilde{S})]$ is generated by the $\mathcal{I}(\tilde{S}' )$ orbits of $\bigcup_{1 \leq i \leq g} [\pi_1(Y_i^1),\pi_1(Y_i^1)]$ \cite[Lemma 4.2]{PM2}. Also, the disk pushing map takes $\bigcup_{1 \leq i \leq g} [\pi_1(Y_i^1),\pi_1(Y_i^1)]$ inside $\bigcup_{1 \leq i \leq g} \mathcal{I}(\Sigma_{g+1,1},Y_i)$. This completes Step 2 and the Lemma.
	\end{proof}
	
	\begin{proof}[Proof of Lemma \ref{cut_handles}]
		If we prove that
		\[ \bigcup_{1 \leq i \leq g} \mathcal{I}(\Sigma_{g+1,1},Y_i) < \bigcup_{I} \mathcal{I}(\Sigma_{g+1,1},X_I). \]
		then we are done, since $\bigcup_{1 \leq i \leq g} \mathcal{I}(\Sigma_{g+1,1},Y_i )$ generates $\mathcal{I}(\Sigma_{g+1,1},S)$ by Lemma \ref{handle_torelli}. We apply induction on genus $g$. The base case is $g=3$ where the Lemma is true. For $g>3$ consider the Birman exact sequence:
		
		\[ 1 \to \pi_1(U S'') \to \mathcal{I}(\Sigma_{g+1,1},Y_i) \to \mathcal{I}(\Sigma_{g+1,1},S'') \to 1. \]
		
		where $S''$ is a surface of genus $g-1$ with 2 boundary components obtained from $Y_i$ after capping $\partial R_i$. We have that the image of $\pi_1(U S'')$ in $\mathcal{I}(\Sigma_{g+1,1},Y_i)$ is inside $\bigcup_I \mathcal{I}(\Sigma_{g+1,1},X_I)$. Also, for $I \subset \boldsymbol{g}$ such that $i \notin I$ we have that $\bigcup_I \mathcal{I}(\Sigma_{g+1,1},X_I)$ is mapped to a generating set in $\mathcal{I}(S'')$ by inductive hypothesis. Therefore, 
		\[ \mathcal{I}(\Sigma_{g+1,1},Y_i) < \bigcup_I \mathcal{I}(\Sigma_{g+1,1},X_I). \]
	\end{proof}
	
	Now we are ready to prove Theorem \ref{cubic}.
	
	\begin{proof}[Proof of Theorem \ref{cubic}]
		Let $\Gamma$ be the group generated by the elements of the Theorem. We proceed to the proof of the Theorem. Let $X_I$ be as in Lemma \ref{cut_handles}. If we glue disks in the boundaries $\partial R_m$ of $X_I$ we get a subsurface namely, $Z$. We note that $Z$ is homeomorphic to a surface of genus $|I|$ with two boundary components $\partial_1$, $\partial_2$. Denote by $\overline{Z}$ the surface obtained by $S$ after gluing a pair of pants in the boundary components $\partial_1$, $\partial_2$. The surface $\overline{Z}$ has genus $|I|+1$ and one boundary component. Note that the map $\I(\Sigma_{ g+1,1},X_I) \to \I(\overline{Z},Z)$ is a well defined epimorphism \cite{PM1}. Consider the Birman exact sequence:
		
		\begin{equation}
			1 \to K \to \I(\Sigma_{ g+1,1},X_I) \to \I(\overline{Z},Z) \to 1.
			\label{ses_stab}
		\end{equation}
		
		We will explain why the generators of $K$ lie in $\Gamma$. For starters, let $Z^i$ be the surface obtained by $X_I$ after capping $\partial R_i$ where $i \notin I$. Denote also $\overline{Z}^i$ the surface obtained from $Z^i$ after gluing a pair of pants to the boundaries $\partial_1, \partial_2$. Hence, we have
		
		\[ 1 \to \pi_1(U Z^i) \to \I(\Sigma_{ g+1,1},X_I) \to \I(\overline{Z}^i,Z^i) \to 1, \]
		
		where $U Z^i$ is the unit tangent bundle of $Z^i$. Note that $i \notin I$, then $\pi_1(U Z^i) < K$
		\[ \big\langle \bigcup_{i \notin I} \pi_1(U Z^i) \big\rangle = K. \]
		
		We identify $\pi_1(U Z^i)$ as the disk pushing map subgroup of $\I(\Sigma_{ g+1,1},X_I)$. We have that $T_{\partial R_i} \in \pi_1(U Z^i)$ and also $\langle T_{\partial R_i}\rangle$ is normal in $\pi_1(U Z^i)$. Also $T_{\partial R_i} \in \Gamma$. If we prove that the generators of $\pi_1(Z^i) = \pi_1(U Z^i) / \langle T_{\partial R_i}\rangle$ we are done. In Figure \ref{BES_kernel_Z} we can see an example of $X_I$ and $Z^i$. \\

		Our aim is to describe the generators of $\pi_1(Z^i)$ as elements inside $\I(\Sigma_{ g+1,1},X_I)$. Recall the arcs $\delta_j$ connecting $p_0$ with $\partial R_j$ and recall also $*_j = \delta_j \cap \partial R_j$. We consider an orientation of $\delta_j$ from $p_0$ to $*_j$. To each $R_j$ consider the two curves $\tilde{x}_j$, $\tilde{y}_j$ as in Figure \ref{disk_pushing_map}. Homotope $\tilde{x}_j$, $\tilde{y}_j$ so that they intersect $*_j$. We define the following curves:
		
		\begin{figure}[h]
			\begin{center}
				\includegraphics[scale=.8]{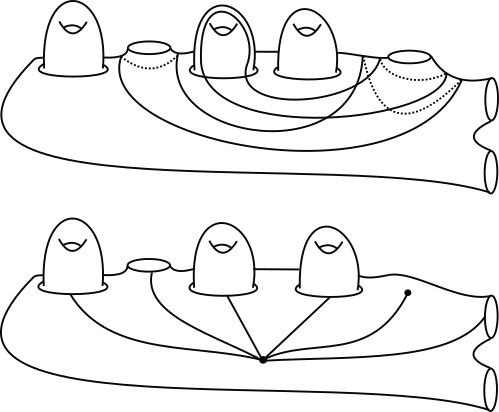}
			\end{center}
			\caption{On the top it is $X_{ \{1,3,4\} }$ and on the bottom it is $Z^5$. Also, on the top we have examples of $\tilde{x}'_3$ and $r_j$.} 
			\begin{picture}(22,12)
				\put(113,125){$*_5$}
				\put(-10,260){$\tilde{x}'_3$}
				\put(-40,240){$r_2$}
			\end{picture}
			\label{BES_kernel_Z}
		\end{figure}
		
		\begin{itemize}
			\item If $j \in I$: $\tilde{x}'_j = \delta_i^{-1} \delta_j \tilde{x}'_j \delta_j^{-1} \delta_i$,
			\item If $j \in I$: $\tilde{y}'_j = \delta_i^{-1} \delta_j \tilde{y}'_j \delta_j^{-1} \delta_i$,
			\item If $j \notin I$: $r_j = \tilde{x}'_j = \delta_i^{-1} \delta_j \partial R_j \delta_j^{-1} \delta_i$,
			\item $r_{g+1} = \tilde{x}'_j = \delta_i^{-1} \delta_j \partial_1 \delta_j^{-1} \delta_i$.
		\end{itemize}
		
		Now replace $*_j$ with $\partial R_j$ to consider the curves above as curves in $X_I$ (see top of Figure \ref{BES_kernel_Z}). Then $\pi_1(Z^i)$ is generated by $T_{\tilde{x}_j} T_{\tilde{x}'_j}^{-1}$, $T_{\tilde{y}_j} T_{\tilde{y}'_j}^{-1}$, $T_{\partial R_j} T_{r_j}^{-1}$, $T_{\partial_1} T_{r_{g+1}}^{-1}$. These elements lie in $\Gamma$ (see also \cite[Proof of Proposition 4.5]{PM3} for similar description of $K$).\\
		
		Note that $\I(\Sigma_{ g+1,1},S_I) < \I(\Sigma_{ g+1,1},X_I)$ since $S_I \subset X_I$. Furthermore, $Z$ is homeomorphic to $S_I$. Recall that one of the boundaries of $S_I$ is $\partial_1$. Denote by $\theta_I$ the other boundary. We homotope the boundary $\partial_2$ of $Z$ to match $\theta_I$. This gives a homotopy from $S_I$ to $Z$, hence $\I(\Sigma_{ g+1,1},S_I) \cong \I(\overline{Z},Z)$. Then, the short exact sequence \ref{ses_stab} is split and 
		\[ \I(\Sigma_{ g+1,1},X_I) \cong \I(\Sigma_{ g+1,1},S_I) \rtimes K \]

	\end{proof}
	
	If we combine Corollary \ref{numerical_thmA} with Theorem \ref{cubic} we deduce Theorem \ref{main_cubic}.
	
	\section{Proof of Theorem C}
	
	Recall that $S$ is a subsurface of $\Sigma_{ g+1,1}$ where $\Sigma_{ g+1,1} \setminus S$ is a pair of pants (recall from Figure \ref{surface_model}). Recall that $\K(\Sigma_{ g+1,1},S) = \I(\Sigma_{ g+1,1},S)\cap \K(\Sigma_{ g+1,1})$; we have seen that $\K(\Sigma_{ g+1,1},S)$ is generated by Dehn twists about $Q$-separating curves (as introduced in Section 2). We want to prove that $\K(\Sigma_{ g+1,1},S)$ is finitely generated. We will use a Theorem by Church-Ershov-Putman \cite[Theorem 3.7]{CEP}.
	
	\begin{theorem}
		Let $G$ be a finitely generated group. Suppose that a group $\Gamma$ acts on $G$ by automorphisms such that the following hold.
		\begin{itemize}
			\item[1.] The group $G$ is generated by a single $\Gamma$-orbit $C\subset G$.
			\item[2.] The image of $\Gamma$ in $\mathrm{GL}(\mathrm{Hom}(G,\mathbb{R}))$ is irreducible in the Zariski topology.
			\item[3.] The graph whose vertices are elements $c \in C$ where two elements are connected by an edge if they commute, is connected.
		\end{itemize}
		Then, $[G,G]$ is finitely generated.
		\label{CEP_THM}
	\end{theorem}
	
	We will use Theorem \ref{CEP_THM} to prove that $[\I(\Sigma_{ g+1,1},S),\I(\Sigma_{ g+1,1},S)]$ is finitely generated. Since $[\I(\Sigma_{ g+1,1},S),\I(\Sigma_{ g+1,1},S)]$ is a subgroup of $\K(\Sigma_{ g+1,1},S)$ of finite index, then we deduce that $\K(\Sigma_{ g+1,1},S)$ is finitely generated. Theorem \ref{CEP_THM} will be used as follows: $\Gamma = \Mod(S)$ and $G = \I(\Sigma_{ g+1,1},S)$.
	
	\begin{definition}
		Denote by BP1 the graph whose vertices are genus 1 $Q$-bounding pair maps in $\I(\Sigma_{ g+1,1},S)$. Two vertices are connected by an edge if the corresponding bounding pair maps commute.
	\end{definition}
	
	\begin{lemma}
		Let $S$ be a surface of genus at least 4 with 2 boundary components. The corresponding graph BP1 is connected if the genus of $S$ is at least 4.
		\label{BP1_connected}
	\end{lemma}
	
	The proof of Lemma \ref{BP1_connected} is the same as in \cite[Proposition 3.3]{CEP} after modifying two things that do not change the hypothesis of the Lemma: Replace the one boundary component with 2 boundary components. Choose the generating set of $\Mod(S)$ depicted in Figure \ref{mcg_gen}. We note that the condition of the genus of $S$ to be at least 4 is necessary. By the change of coordinates principle (\cite[Chapter 1]{BFM}) that for any two Q-bounding pair maps we can find an orientation preserving homeomorphism of $S$ that sends one Q-bounding pair map to the other. If we combine the latter fact with Lemmas \ref{bp_maps} and \ref{bp_genus1} we deduce that $\I(\Sigma_{ g+1,1},S)$ is generated by a single $\Mod(S)$ orbit of an element of BP1. We have proved the following:
	
	\begin{itemize}
		\item[1.] The group $\I(\Sigma_{ g+1,1},S)$ (which is finitely generated) is generated by a single $\Mod(S)$-orbit of a Q-bounding pair map of $\I(\Sigma_{ g+1,1},S)$.
		\item[3.] The graph of BP1 is connected.
	\end{itemize}
	
	It remains to prove that the image of $\Mod(S)$ in $\mathrm{GL}(\mathrm{Hom}(\I(\Sigma_{ g+1,1},S),\mathbb{R}))$ is irreducible in the Zariski topology.
	
	\subsection{Zariski topology}

	We will define a map $\Mod(S) \to GL\big( \mathrm{Hom}(\I(\Sigma_{ g+1,1},S),\Rr) \big)$. The Zariski topology of $\Mod(S)$ is the pullback Zariski topology of $GL\big( \mathrm{Hom}(\I(\Sigma_{ g+1,1},S),\Rr) \big)$. For our arguments below we follow Humphrey's book \cite{HU1975}.
	
	\begin{itemize}
		\item[Fact 1:] Let $Y$ be irreducible, and let $f: X\to Y$ be any function. Give $X$ the pullback topology (closed sets in $X$ are exactly preimages of closed sets in $Y$). Then $X$ is irreducible.
		
		\item[Fact 2:] If $X$ is irreducible and $f: X \to Y$ is continuous, then $f(X)$ (with the subspace topology from $Y$) is irreducible.
		
		\item[Fact 3:] Let $Z$ be a subspace of a topological space $W$, and let $\overline{Z}$ be its closure in $W$. Then $Z$ is irreducible if and only if $\overline{Z}$ is irreducible.
		
		\item[Fact 4:] An algebraic group is irreducible if and only if it is connected.
		
		\item[Fact 5:] Two algebraic groups $X,Y$ are irreducible (or connected by Fact 4) then $X \times Y$ is irreducible (or connected).
	\end{itemize}
	
	Now we turn into the mapping class group. As explained in Section 2 we have a surjective homomorphism
	\[ \Mod(S) \to (\mathrm{Sp}_{2g+2}(\Zz))_{y_{g+1}}, \]
	where $(\mathrm{Sp}_{2g+2}(\Zz))_{y_{g+1}}$ is the stabilizer of $y_{g+1}$ in $\mathrm{Sp}_{2g+2}(\Zz)$ (Note that it could be the stabilizer of any vector in $\mathrm{H}_1(S,Q)$; $y_{g+1}$ is a convenient choice). We have that $\mathrm{Sp}_{2g+2}(\Zz) \subset \mathrm{Sp}_{2g+2}(\Rr)$. \\
	
	Recall from Section 2 that
	\[V_\Zz = \Zz x_1 \oplus \cdots \oplus \Zz x_g \bigoplus \Zz y_1 \oplus \cdots \oplus \Zz y_g \oplus \Zz y_{g+1}\]
	where $x_i, y_i$ are the basis vectors of $\mathrm{H}_1(S,Q)$ (except $x_{g+1}$). If we replace $\Zz$ by $\Rr$ we obtain $V_\Rr$.
%
	Recall the Zariski continuous maps from \cite[Section 3]{CEP}
	\[ \Mod(\Sigma_{ g+1,1}) \to \Sp{2g+2}{\Zz} \to \Sp{2g+2}{\Rr} \to GL\big( \bigwedge^3 \mathrm{H}_1(\Sigma_{ g+1,1};\Rr) \big) \to GL\big( \mathrm{Hom}(\bigwedge^3 \mathrm{H}_1(\Sigma_{ g+1,1};\Rr),\Rr) \big). \]
	Since $\Mod(S) < \Mod(\Sigma_{ g+1,1})$ and $(\Sp{2g+2}{R})_{y_{g+1}} < \Sp{2g+2}{R}$ where $R=\Zz, \Rr$, then the sequence of maps above gives a sequences of continuous maps below
	\begin{equation}
		\Mod(S) \to (\Sp{2g+2}{\Zz})_{y_{g+1}} \to (\Sp{2g+2}{\Rr})_{y_{g+1}} \to GL\big( \bigwedge^3 V_{\Rr} \big) \to GL\big( \mathrm{Hom}(\bigwedge^3 V_{\Rr},\Rr) \big).
		\label{cont_maps}
	\end{equation}
	In order to complete the proof of Theorem C we only need to prove the following Lemma.
	
	\begin{lemma}
		The image of $\Mod(S)$ in $GL\big( \mathrm{Hom}(\bigwedge^3 V_{\Rr},\Rr) \big)$ is Zariski irreducible.
	\end{lemma}
	
	\begin{proof}
		The homomorphism $\Mod(S) \to (\Sp{2g+2}{\Zz})_{y_{g+1}}$ is surjective. The closure of $(\Sp{2g+2}{\Zz})_{y_{g+1}}$ is $(\Sp{2g+2}{\Rr})_{y_{g+1}}$. Hence, if $(\Sp{2g+2}{\Rr})_{y_{g+1}}$ is irreducible, then by Fact 3, then, $(\Sp{2g+2}{\Zz})_{y_{g+1}}$ is irreducible. Consequently, by Fact 1, $\Mod(S)$ is irreducible in the Zariski topology. Hence, we need to prove that $(\Sp{2g+2}{\Rr})_{y_{g+1}}$ is irreducible. By Section 2, we have that $(\Sp{2g+2}{\Rr})_{y_{g+1}}$ as a set is represented by $\Sp{2g+2}{\Rr} \times \Rr ^{2g+1} \times \Rr$. Each of the factors is connected. The Cartesian product of connected spaces is connected by Fact 5. Thus, $(\Sp{2g+2}{\Rr})_{y_{g+1}}$ is connected and hence irreducible by Fact 4. Since all maps of (\ref{cont_maps}) are continuous, then we are done by Fact 2.
	\end{proof}
	
	\bibliographystyle{plain}
	\bibliography{bibliography.bib}{}
\end{document}